\documentclass[a4paper,11pt]{amsart}

\usepackage{amsmath}
\usepackage{amstext}
\usepackage{amsbsy}
\usepackage{amsopn}
\usepackage{amsthm}
\usepackage{amscd}
\usepackage{amsxtra}
\usepackage{amsfonts}
\usepackage{pictexwd}
\usepackage{ifthen}
\usepackage{xspace}
\usepackage{fancyheadings}
\usepackage{layout}
\usepackage{hyperref}

\input xy
\xyoption{matrix}
\xyoption{arrow}

\newdir{ (}{{}*!/-5pt/@^{(}}

\newcommand{\ZZ}{\mathbb{Z}}

\newcommand{\PZ}{\mathbb{P}}

\newcommand{\GZ}{\mathbb{G}}

\newcommand{\SZ}{\mathbb{S}}

\newcommand{\sO}{{\mathcal O}}
\newcommand{\sL}{{\mathcal L}}

\newcommand{\sI}{{\mathcal I}}
\newcommand{\sT}{{\mathcal T}}

\newcommand{\from}{\leftarrow}

\newboolean{probleme}
\newcommand{\problem}[1]
           {\ifthenelse{\boolean{probleme}}
                       {{\bf(PROBLEM: #1)\bf}}
                       {}
           }
\newboolean{zukuenftiges}
\newcommand{\zukunft}[1]
           {\ifthenelse{\boolean{zukuenftiges}}
                       {{\bf(AUSBAUM\"OGLICHKEIT: #1)\bf}}
                       {}
           }
\newboolean{extras}
\newcommand{\extra}[1]
           {\ifthenelse{\boolean{extras}}
                       {{\bf EXTRA #1 EXTRA\bf}}
                       {}
           }

\DeclareMathOperator{\codim}{codim}

\DeclareMathOperator{\Img}{Im}

\DeclareMathOperator{\Hom}{Hom}

\DeclareMathOperator{\rank}{rank}

\DeclareMathOperator{\cliff}{cliff}
\DeclareMathOperator{\Syz}{Syz}
\DeclareMathOperator{\Gensyz}{Gensyz}

\DeclareMathOperator{\id}{id}
\DeclareMathOperator{\GL}{GL}

\DeclareMathOperator{\Flag}{Flag}

\DeclareMathOperator{\Gr}{Gr}                 
\DeclareMathOperator{\pfaff}{pfaff}

\newcommand{\youngbox}[2]{\mspace{-0.7mu}\framebox[#1]{\rule{0mm}{1.5ex}$#2$}}

\newcommand{\xyoungbox}[1]{\mspace{-0.7mu}\framebox{\rule{0mm}{1.5ex}$#1$}}
\newcommand{\blackbox}{\mspace{-0.7mu}\rule[-3.4pt]{5mm}{13.9pt}}
\newcommand{\ybox}[1]{\youngbox{5mm}{#1}}

\swapnumbers
\theoremstyle{plain}
\begingroup

\newtheorem{thm}{Theorem}
\newtheorem{cor}[thm]{Corollary}
\newtheorem{lem}[thm]{Lemma}
\newtheorem{lemdefn}[thm]{Lemma and Definition}
\newtheorem{propdefn}[thm]{Proposition and Definition}
\newtheorem{prop}[thm]{Proposition}
\newtheorem{conj}[thm]{Conjecture}

\numberwithin{thm}{subsection} 
\endgroup

\theoremstyle{definition}
\newtheorem{defn}[thm]{Definition}
\newtheorem{rem}[thm]{Remark}

\numberwithin{equation}{section}

\newcommand{\nosubsections}{\renewcommand{\thethm}{\thesection.\arabic{thm}}
                            \setcounter{thm}{0}
                           }

\newcommand{\cref}[3]{(\ref{#1}, #2 \ref{#3})}

\date{\today}
\email{bothmer@btm8x5.mat.uni-bayreuth.de}
\address{Hans-Christian Graf v. Bothmer\\ Mathematisches Institut der Universit\"at
Bayreuth\\ 95440 Bayreuth\\ Germany }

\setboolean{probleme}{true}

\begin{document}

\title[Geometric Syzygies of Canonical Curves]{Geometric Syzygies of Canonical Curves of even Genus lying on
       a $K3$-Surface}

\thanks{Supported by the Schwerpunktprogramm ``Global Methods in Complex
        Geometry'' of the Deutsche Forschungs Gemeinschaft}

\author{Hans-Christian Graf v. Bothmer}
\maketitle

\setcounter{tocdepth}{1}
\tableofcontents

\section{Introduction}
\nosubsections

\newcommand{\PZorth}{\PZ^\bot}

In this paper we study the syzygies of canonical curves
$C \subset \PZ^{g-1}$ for $g = 2k$.

In \cite{GL1} Green and Lazarsfeld construct low-rank-syzygies of $C$ from 
special linear systems on $C$. More precisely linear systems
of Clifford index $c$ give a $(g-c-3)$-rd syzygy. We call these syzygies
geometric syzygies. Green's conjecture
\cite{GreenKoszul} paraphrased in this way is

\begin{conj}[Green]
Let $C$ be a canonical curve, then 
\[
\text{$C$ has no geometric $p$-th syzygies}
\iff
\text{$C$ has no $p$-th syzygies at all}
\]
\end{conj}

This conjecture as received a lot of attention in the last years
and it is now known in many cases \cite{Petri}, \cite{GreenKoszul}, \cite{Sch86},
\cite{Sch88}, \cite{Voisin}, \cite{Sch91}, \cite{Ehbauer}, \cite{HirschR},
\cite{TeixidorGreen}, \cite{VoisinK3}.

Much less is know about the following 
natural generalization of Green's conjecture:

\begin{conj}[Geometric Syzygy Conjecture]
Let $C$ be a canonical curve, then the scheme of geometric $p$-th 
syzygies spans the
space of all $p$-th syzygies.
\end{conj}

For general canonical curves both conjectures 
are equivalent in the range $p \ge \frac{g-3}{2}$ since a 
general canonical curve has no linear systems of Clifford index
$c \le \frac{g-3}{2}$.

The case $p=0$ (geometric quadrics) of the geometric syzygy conjecture
was proved by \cite{AM} for general canonical curves, 
and by \cite{GreenQuadrics}
for all canonical curves. The case $p=1$ was done for general canonical curves
of genus $g \ge 9$ in \cite{HC}.

In this paper, based on Voisin's recent result \cite{VoisinK3},
we clarify the cases with $g=2k$ and $p=k-2$ for
canonical curves lying on general $K3$ surfaces.

\begin{figure}[h]
\mbox{
\beginpicture   
\setlinear 

\setcoordinatesystem units <6mm,6mm> point at -150 0.5
\unitlength6mm
\setplotarea x from 0 to 16, y from 0 to 9
\axis left 
   ticks numbered from 0 to 9 by 1 /
\put{\vector(0,1){10}} [Bl] at 0 0
\put{$p$} at 0 10.5
\axis bottom
   ticks numbered from 0 to 15 by 1 /
\put {\vector(1,0){16}} [Bl] at 0 0 
\put{$g$} at 16.5 0

\put{\line(1,1){9}} [Bl] at 4 0
\put{\line(2,1){11}} [Bl] at 4 0

\put{\circle*{0.1}} [Bl] at 0.7 9
\put{\lines{Green's Conj. $\implies$ GSC}} [l] at 1.2 9

\setplotsymbol ({\circle{0.05}} [Bl])
\setsolid
\plot 0.6 8  0.7 8    0.8 8 /
\put{\lines{known by \cite{AM}, \cite{GreenQuadrics}, \cite{HC}}} [l] at 1.2 8

\putrectangle corners at 0.6 6.9 and 0.8 7.1
\put{\lines{known by \cite{HCandRanestad}}}[l] at 1.2 7

\put{\circle{0.2}} [Bl] at 0.7 6
\put{\lines{unknown}} [l] at 1.2 6

\put{\lines{this article}} [l] at 1.2 5

\setplotsymbol ({\circle{0.05}} [Bl])
\setsolid
\plot 4 0  8 0  15 0 /

\plot 9 1  10 1  15 1 /

\putrectangle corners at  5.9 0.9 and  6.1 1.1
\putrectangle corners at  7.9 1.9 and  8.1 2.1
\putrectangle corners at  6.9 0.9 and  7.1 1.1

\linethickness1.2mm 
\putrule from  9.9 3 to  10.1 3
\putrule from 11.9 4 to  12.1 4
\putrule from 13.9 5 to  14.1 5
\putrule from  0.6 5 to  0.8 5

\put{\circle{0.2}} [Bl] at 8 1
\put{\circle{0.2}} [Bl] at 9 2
\put{\circle{0.2}} [Bl] at 10 2
\put{\circle{0.2}} [Bl] at 11 2
\put{\circle{0.2}} [Bl] at 12 2
\put{\circle{0.2}} [Bl] at 13 2
\put{\circle{0.2}} [Bl] at 14 2
\put{\circle{0.2}} [Bl] at 15 2
\put{\circle{0.2}} [Bl] at 11 3
\put{\circle{0.2}} [Bl] at 12 3
\put{\circle{0.2}} [Bl] at 13 3
\put{\circle{0.2}} [Bl] at 14 3
\put{\circle{0.2}} [Bl] at 15 3
\put{\circle{0.2}} [Bl] at 13 4
\put{\circle{0.2}} [Bl] at 14 4
\put{\circle{0.2}} [Bl] at 15 4
\put{\circle{0.2}} [Bl] at 15 5

\put{\circle*{0.1}} [Bl] at 5 1
\put{\circle*{0.1}} [Bl] at 6 2
\put{\circle*{0.1}} [Bl] at 7 2
\put{\circle*{0.1}} [Bl] at 7 3
\put{\circle*{0.1}} [Bl] at 8 3
\put{\circle*{0.1}} [Bl] at 9 3
\put{\circle*{0.1}} [Bl] at 8 4
\put{\circle*{0.1}} [Bl] at 9 4
\put{\circle*{0.1}} [Bl] at 10 4
\put{\circle*{0.1}} [Bl] at 11 4
\put{\circle*{0.1}} [Bl] at 9 5
\put{\circle*{0.1}} [Bl] at 10 5
\put{\circle*{0.1}} [Bl] at 11 5
\put{\circle*{0.1}} [Bl] at 12 5
\put{\circle*{0.1}} [Bl] at 13 5
\put{\circle*{0.1}} [Bl] at 10 6
\put{\circle*{0.1}} [Bl] at 11 6
\put{\circle*{0.1}} [Bl] at 12 6
\put{\circle*{0.1}} [Bl] at 13 6
\put{\circle*{0.1}} [Bl] at 14 6
\put{\circle*{0.1}} [Bl] at 15 6
\put{\circle*{0.1}} [Bl] at 11 7
\put{\circle*{0.1}} [Bl] at 12 7
\put{\circle*{0.1}} [Bl] at 13 7
\put{\circle*{0.1}} [Bl] at 14 7
\put{\circle*{0.1}} [Bl] at 15 7
\put{\circle*{0.1}} [Bl] at 12 8
\put{\circle*{0.1}} [Bl] at 13 8
\put{\circle*{0.1}} [Bl] at 14 8
\put{\circle*{0.1}} [Bl] at 15 8
\put{\circle*{0.1}} [Bl] at 13 9
\put{\circle*{0.1}} [Bl] at 14 9
\put{\circle*{0.1}} [Bl] at 15 9

\endpicture   
     }
   \caption{The geometric syzygy conjecture for general curves.} 
\end{figure}

We proceed as follows:

By a construction of Mukai \cite{ClassificationMukai} 
a general $K3$-surface $S$ of sectional genus
$g = 2k$ can be embedded in a Grassmannian $\GZ = \Gr(k+2,2)$, such that
\[
         S \subset \GZ \cap \PZ^g
\]
for a particular $\PZ^g \not\subset \GZ$. We reconstruct this embedding
from a $(k-2)$-nd Grassmannian syzygy of $S$.

The minimal free resolution
of $\GZ$ is known by work of J\'ozefiak, Pragacz and Weyman
\cite{PfaffianResolution}. It has the form:

\begin{center}
\mbox{
\beginpicture   
\setlinear 

\setcoordinatesystem units <5mm,5mm> point at -150 0.5
\unitlength5mm
\setplotarea x from 0 to 22, y from 0 to 7

\putrectangle corners at   0 7 and   1 6
\putrectangle corners at   1 6 and   6 5
\putrectangle corners at   6 6 and   7 5
\putrectangle corners at   3 5 and  12 4
\put{$\ddots$}  at 11 3.7
\putrectangle corners at  10 3 and  19 2
\putrectangle corners at  15 2 and  21 1
\putrectangle corners at  21 1 and  22 0

\put{$*$} at 6.5 5.5

\put{$\overbrace{\rule{30mm}{0mm}}^{k-1}$}[Bl] at 1 6

\endpicture   
     }
\end{center}

using the MACAULAY-notation \cite{M2}. The minimal free resolution of the $K3$
surface is also known by a recent result of Voisin \cite{VoisinK3}:

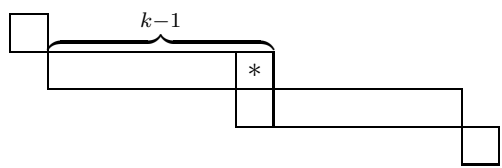
\begin{center}
\mbox{
\beginpicture   
\setlinear 

\setcoordinatesystem units <5mm,5mm> point at -150 0.5
\unitlength5mm
\setplotarea x from 0 to 22, y from 3 to 7

\putrectangle corners at   0 7 and   1 6
\putrectangle corners at   1 6 and   6 5
\putrectangle corners at   6 6 and   7 5
\putrectangle corners at   6 5 and   7 4
\putrectangle corners at   7 5 and  12 4
\putrectangle corners at  12 4 and  13 3

\put{$*$} at 6.5 5.5

\put{$\overbrace{\rule{30mm}{0mm}}^{k-1}$}[Bl] at 1 6

\endpicture   
     }
\end{center}

Our main observation is that both resolutions have a linear strand
of the same length and that the dimensions of the last nonzero linear
syzygy spaces $(*)$ are equal.

Using Voisin's theorem we show that the natural map between the
two spaces is an isomorphism even though the intersection $\GZ \cap \PZ^g$
is not of expected dimension. It turns out to be enough, that $S$ is
irreducible.

We go on to describe the space of minimal rank
syzygies $Y_{min}$ of $\GZ$ and $S$ as a $(k-2)$-uple embedded
$\PZ^{k+1}$.

Cutting down one more dimension to a canonical curve $C \subset \PZ^{g-1}$
we obtain a finite number of lines of geometric syzygies
in $\PZ^{k+1}$ whose image under the $(k-2)$-uple embedding spans the
space of all $(k-2)$-nd syzygies of $C$.

As a main tool of our work we associate two geometric objects to a
syzygy $s$:
\begin{itemize}
\item[(1)] The space of linear forms $L_s$ involved in $s$. This 
space and its associated vector bundle of linear forms allows us to control
the rank of $s$ after we have cut down to $S$ or $C$.  
\item[(2)] The syzygy scheme $\Syz(s)$ of $s$ which is in a certain sense
the vanishing locus of $s$. It is used to prove that certain
syzygies do survive the restriction to the linear subspaces $\PZ^g$ 
and $\PZ^{g-1}$.
\end{itemize}

Several of our arguments also appear in Voisin's proof of her theorem,
but with different aims. Most notably she also shows that
$\PZorth \cap \GZ^* = \emptyset$ and that the above rational normal
curves span the space of the $(k-2)$uple embedding. As far as we know
our corollary \ref{injective} is new. 

I would to thank Frank-Olaf Schreyer for introducing me to this
subject and for the many helpful discussions leading to this paper. 
Also I am grateful to Thomas Eckl, for the numerous discussions  
which clarified many details of this work.

\section{Syzygies of Low Rank}
\nosubsections

Let $X \subset \PZ^n \cong \PZ(V)$ be a irreducible non degenerate variety,
$I_X$ generated by quadrics and
\[
   \sI_X \from V_0 \otimes \sO(-2) 
         \xleftarrow{\varphi_1} V_1 \otimes \sO(-3)
         \xleftarrow{\varphi_2} \dots
         \xleftarrow{\varphi_m} V_m \otimes \sO(-m-2)
\]
the linear strand of its minimal free resolution

\begin{defn}
An element $s \in V_p$ is called a $p$-th (linear) syzygy of $X$.
$\PZ(V_p^*)$ is called the space of $p$-th syzygies.
\end{defn}

Every linear syzygy $s$ involves a well defined number of linearly independent
linear forms. This number is called the rank of $s$. In a more formal way we 
have:

\newcommand{\tildephi}{\tilde{\varphi}}
\newcommand{\spans}{L_s}

\begin{defn}\label{tildephi}
Let $s \in V_p$ be a syzygy, 
\[
   \tildephi \colon V_p \to V_{p-1}\otimes V
\]
the map of vector spaces induced by $\varphi_p$. Then the image
of $s$ under $\tildephi$ can be interpreted as a linear
map:
\[
   \tildephi(s) \in V_{p-1} \otimes V \cong \Hom(V_{p-1}^*,V).
\]
With this
\[
    \spans := \Img \tildephi(s) \subset V
\]
is called the space of linear forms involved in $s$, and
\[
    \rank s := \rank \tildephi(s) = \dim \spans
\]
the rank of $s$.
\end{defn}

To apply geometric methods to the study of low rank syzygies we projectivize
the space of $p$-th syzygies to $\PZ(V_p^*)$ and 
give a determinantal description
of the space $Y_{min}$ of minimal rank syzygies. 
The linear forms involved in these
syzygies define a vector bundle on $Y_{min}$:

\begin{defn}\label{mapofvectorbundles}
On the space of $p$-th syzygies $\PZ(V_p^*)$ the map of
vector spaces $\tildephi_p$ induces 
a map of vector bundles
\[
    \psi \colon V_{p-1}^* \otimes \sO_{\PZ(V_p^*)} (-1)
                \to V \otimes  \sO_{\PZ(V_p^*)}
\]
that satisfies
\[
       \psi|_s = \tildephi_p(s) \in \Hom(V_{p-1}^*,V)
\]
The determinantal loci $Y_r(\psi) \subset \PZ(V_p^*)$ of $\psi$ are called
schemes of rank $r$ syzygies, since the syzygies in their support
have rank $\le r$.

On the scheme of minimal rank syzygies $Y_{min} := Y_{r_{min}}(\psi)$ 
the restricted
map $\psi|_{Y_{min}}$ has constant rank $r_{min}$. Therefore the image
$\sL := \Img(\psi|_{Y_{min}})$ is a vector bundle. We call it the vector bundle
of linear forms, since
\[
                 \sL|_s = \spans \subset V 
\]
for all minimal rank syzygies $s \in Y_{min}$.
\end{defn}

\section{Syzygy Schemes}
\nosubsections

A second geometric object associated to a syzygy $s$ is
obtained by calculating $V_p$ via Koszul cohomology:

\begin{lem} \label{koszul}
\[
V_p \cong  \ker \bigl( \Lambda^p V \otimes (I_X)_2 \to 
                        \Lambda^{p-1} V \otimes (I_X)_3 
                 \bigr)
     \cong  H^0 \bigl( \PZ(V),\Omega^p(p+2)\otimes I_X \bigr)
\]
\end{lem}

\begin{proof} \cite{GreenKoszul}, \cite{Ehbauer}
\end{proof}

So linear syzygies are twisted $p$-forms that vanish on $X$. Often 
these $p$-forms vanish on a larger variety:

\begin{defn}
Let $s \in V_p$ be
a $p$-th syzygy of $X$. Then the vanishing set $\Syz(s)$ of the
corresponding twisted $p$-form is called the syzygy scheme of $s$.
\end{defn}

The ideal of a syzygy scheme can be calculated via:

\begin{lemdefn}\label{quadricsinvolved}
Let $\{v_\alpha\}$ be a basis of $\Lambda^p V$. Then every syzygy
$s \in V_p$ can be uniquely written as
\[
           s = \sum_\alpha v_\alpha \otimes Q_\alpha
\]
where $Q_\alpha$ are quadrics in the ideal of $X$, and
the ideal of $\Syz(s)$ is generated by the $Q_\alpha$. This 
ideal is also called the ideal of quadrics involved in $s$.
\end{lemdefn}

\begin{proof}
Since $V_q = \ker \bigl( \Lambda^p V \otimes (I_X)_2 \to 
                       \Lambda^{p-1} V \otimes (I_X)_3 
                 \bigr)$
every syzygy can be written as above. Since the $v_\alpha$ are
linearly independent, $s = \sum_\alpha v_\alpha \otimes Q_\alpha$ vanishes
if and only if all $Q_\alpha$ vanish.
\end{proof}

Often the syzygies of low rank have the most interesting
syzygy varieties. Some of them can be calculated with the methods
of the next section.

\section{Generic Syzygy Schemes}
\nosubsections

We now consider syzygies $s$ of low rank $r$ and their syzygy schemes. 
At it turns out, these syzygy schemes are always cones over
linear sections of certain generic syzygy schemes:

\newcommand{\genspace}{ L \oplus \Lambda^{r-p-1} L }

\begin{defn}
Let $L$ be an $r$-dimensional vector space. Then
\[
       \Gensyz_p(L) = \{ (l^*,a^*) \in L^* \oplus \Lambda^{r-p-1} L^*
                         \, | \, l^* \wedge a^* = 0 \}
                    \subset \PZ(\genspace)
\]
is called the $p$-th generic syzygy scheme of $L$.
\end{defn}

This definition goes back to Eusen and Schreyer \cite{Eusen}. 
We will now recall some 
well known facts about generic
syzygy schemes.

First of all, the ideal of a generic syzygy variety can be easily 
calculated:

\begin{prop} \label{GenSyzEquations}
$\Gensyz_p(L)$ is cut out by the quadrics in the image of
\[
        \Lambda^{r-p} L \to L\otimes \Lambda^{r-p-1} L 
                        \hookrightarrow S_2(\genspace).
\]
\end{prop}

\begin{proof}
Consider the natural map
\[
\begin{array}{rr@{\otimes}lcc}
    \varphi \colon  
     &L^* &\Lambda^{r-p-1} L^* &\to  &\Lambda^{r-p} L^* \\
     &l^* &a^*               &\mapsto &l^* \wedge a^*
\end{array}
\]
Notice that $(l^*, a^*)$ is in $\Gensyz_p(L)$ if and only if
$(l^* \otimes a^*)$ is in the kernel of $\varphi$. Now $\ker \varphi$
is cut out by the linear forms in the image of the dual map

\[
    \varphi^* \colon  \Lambda^{r-p} L \to  L \otimes  \Lambda^{r-p-1} L 
\]

On $\PZ(L \oplus \Lambda^{r-p-1} L)$ these linear forms
become quadrics that cut out $\Gensyz_p(L)$.            
\end{proof}

We now identify an up to a constant canonically defined 
$p$-th syzygy $s_{gen}$ of $\Gensyz_p(L)$ with $L_{s_{gen}}=L$:

\begin{propdefn} \label{genericsyzygy}
Let $L$ be an $r$-dimensional vector space, and $s_{gen} \in \Lambda^r L$
an orientation. Then $s_{gen}$ is a natural $p$-th syzygy of 
$\Gensyz_p(L)$ via the inclusion
\[
             \Lambda^r L \hookrightarrow 
             \Lambda^p L \otimes (L \otimes \Lambda^{r-p-1} L).
\]
Furthermore
\[
          \Syz(s_{gen}) = \Gensyz_p(L).
\]
We call $s_{gen}$ a generic $p$-th syzygy.
\end{propdefn}

\begin{proof}
First we show, that $s_{gen}$ is really a $p$-th syzygy of $\Gensyz_p(L)$. 
For this we have the following maps

\begin{center}
\mbox{
\xymatrix{
            {\begin{array}{c} 
                { \scriptstyle s \in } \\
                \Lambda^r L
             \end{array}}  
            \ar@<-1ex>[r]^-1 
            \ar[d]^-2
            & {\begin{array}{c} 
                { \scriptstyle s' \in} \\
                \Lambda^p \otimes \Lambda^{r-p} L
             \end{array}}  
            \ar[d]^-3 \\
            \Lambda^{p+1} L \otimes \Lambda^{r-p-1} L \ar[r]^-4
            & 
            {\begin{array}{c} 
                { \scriptstyle s_{gen} \in } \\
                \Lambda^p L \otimes (L \otimes \Lambda^{r-p-1} L) 
             \end{array}}  
            \ar[r]^-5 \ar[d]
            & \Lambda^{p-1} \otimes (S_2(L) \otimes \Lambda^{r-p-1} L) \ar[d] \\
            & \Lambda^p V \otimes S_2(V) \ar[r]^-6
            & \Lambda^{p-1} V \otimes S_3(V)
         }
      }
\end{center}

Mapping $s$ 
via $1$ to $s' $
and then via $3$ to $s_{gen}$
shows
\[ 
           s_{gen} \in \Lambda^{p} L \otimes (I_{\Gensyz_p(L)})_2
\]
since 
\[ 
           (I_{\Gensyz_p(L)})_2 = 
           \Img(\Lambda^{r-p} L \to L \otimes \Lambda^{r-p-1} L)
\]
by proposition \ref{GenSyzEquations}.

Mapping $s$ via $2$ and $4$ to $s_{gen}$ shows that $s_{gen}$ is in the
kernel of $5$ since the middle row of the above diagram is a complex.
Now $5$, which is the restriction of $6$, restricts further to
\[
      \varphi \colon     
           \Lambda^{p} L \otimes (I_{\Gensyz_p(L)})_2 \to
           \Lambda^{p} L \otimes (I_{\Gensyz_p(L)})_3
\]
and $s_{gen}$ is consequently also in the kernel of $\varphi$. This proves
that $s_{gen}$ is a $p$-th syzygy of $\Gensyz_p(L)$.

For $\Syz(s_{gen})=\Gensyz_p(L)$ notice that $s'$ is a trace element
\[
       s'  = \sum w_\alpha^* \otimes w_\alpha 
             \in  (\Lambda^{r-p} L)^* \otimes \Lambda^{r-p} L
\]
with $\{w_\alpha\}$ a basis of $\Lambda^{r-p} L$ and 
$(\Lambda^{r-p} L)^* \cong \Lambda^p L$ via the orientation $s$. Now
since $\{w_\alpha\}$ is a basis of $\Lambda^{r-p} L$, the quadrics involved
in image $s_{gen}$ of $s'$ form a basis of
\[ 
           \Img(\Lambda^{r-p} L \to L \otimes \Lambda^{r-p-1} L)
           = (I_{\Gensyz_p(L)})_2.
\]
This proves
$\Syz(s_{gen})=\Gensyz_p(L)$.
\end{proof}

We are now ready to stated the main result of this section:

\begin{thm} \label{maptogensyz}
Let $X \subset \PZ(V)$ be a non degenerate, possibly reducible variety,
$I_X$ generated by quadrics and $s \in V_p$ a $p$-th syzygy of $X$.
If we denote the space of linear forms involved in $s$ by $L_s$, then
there exists a linear map
\[
         \pi^* \colon L_s \oplus \Lambda^{r-p-1} L_s \to V
\]
and a generic $p$-th syzygy $s_{gen}$ of $\Gensyz_p(L_s)$ such that
\[
            s = \pi^* (s_{gen})
\]
More geometrically consider the linear projection
\[
            \pi \colon \PZ(V) \dasharrow \PZ^n 
                       \subset \PZ(L_s \oplus  \Lambda^{r-p-1} L_s)
\]
corresponding to $\pi^*$. Then
\[
        \Syz(s) = \pi^{-1}(\PZ^n \cap \Gensyz_p(L_s)).
\]
\end{thm}

\begin{proof}
To construct $\pi^*$ we have to make the isomorphism
\[
        V_p \cong \ker(\Lambda^p V \otimes (I_X)_2
                       \to \Lambda^{p-1} V \otimes (I_X)_3)
\]
from lemma \ref{koszul} more explicit. 

Let 
\[
(*)\quad\quad 
         \sO_{\PZ(V)} \from V_0 \otimes \sO(-2)
                      \from \dots
                      \from V_{p-1} \otimes \sO(-p-1)
                      \xleftarrow{\varphi} V_p \otimes \sO(-p-2)
\] 
be the linear strand of the minimal free resolution of $X$. As in
definition \ref{tildephi} $\varphi$ induces a map
\[
   \tildephi \colon V_p \to V_{p-1} \otimes V \cong \Hom(V_{p-1}^*,V)
\]
by taking global sections. $\tildephi(s)$ is a map
from $V_{p-1}^*$ to $V$ with image $L_s$.

$s$ and $\tildephi(s)$ induce a map between the dual of $(*)$ and the
Koszul-complex associated to $L_s$

\begin{center}
\mbox{
\xymatrix{
            \sO \ar[d]^{\alpha} \ar[dr]^\sigma \ar[r]
            & V_0^* \otimes \sO(2) \ar[d]^{\tildephi_p(s)} \ar[r]
            & \dots \ar[r]
            & V_{p-1}^* \otimes \sO(p+1) \ar[d]^{\tildephi(s)} \ar[r]
            & V_{p-2}^* \otimes \sO(p+2) \ar[d]^{s^*}\\
            \Lambda^{p+1} L_s \otimes \sO(1)  \ar[r]
            & \Lambda^{p} L_s \otimes \sO(2) \ar[r]
            & \dots \ar[r]
            & L_s \otimes \sO(p+1) \ar[r]
            & \sO(p+2)
         }
     }
\end{center}

All liftings of $\tildephi(s)$ except $\alpha$ are of degree $0$. Since 
the Koszul-complex is a minimal free resolution, all these liftings are
uniquely determined. Consequently $\sigma$ is also uniquely determined.

$\sigma$ is a section in 
\[
     H^0\bigl(\Lambda^{p} L_s \otimes \sO(2)\bigr) = \Lambda^{p} L_s \otimes S_2(V)
\]
and since it factors over $V_0^* \otimes \sO(2)$ all quadrics involved
in $\sigma$ are in $(I_X)_2$. On the other hand $\sigma$ factors
over $\Lambda^{p+1} L_s \otimes \sO(1)$ and therefore
\[
    \sigma \in \ker\bigl(\Lambda^p L_s \otimes (I_X)_2 \to 
                    \Lambda^{p-1} L_s \otimes (I_X)_3\bigr)
       \subset \ker\bigl(\Lambda^p V \otimes (I_X)_2 \to 
                    \Lambda^{p-1} V \otimes (I_X)_3\bigr).
\]
This defines a map
\[
     V_p \to \ker\bigl(\Lambda^p V \otimes (I_X)_2 \to 
                    \Lambda^{p-1} V \otimes (I_X)_3\bigr).
\]
For the inverse map we take
\[
     \sigma \in  \ker\bigl(\Lambda^p V \otimes (I_X)_2 \to 
                    \Lambda^{p-1} V \otimes (I_X)_3\bigr)
\]
and dualize to
\begin{center}
\mbox{
\xymatrix{
            \sO 
            & V_0 \otimes \sO(-2) \ar[l] \\
            & \Lambda^{p} L_s^* \otimes \sO(-2). \ar[lu]^{\sigma^*} \ar[u] 
         }
     }
\end{center}
Now $\sigma^*$ lifts to a map of complexes from the Koszul complex
associated to $L_s^*$ to the linear strand of $X$. The last map
\[
             \Lambda^0 L_s^* \otimes \sO(-p-2) \to V_p \otimes \sO(-p-2)
\]
produces a unique $p$-th syzygy.

Our projection $\pi$ is now constructed from the section 
$\alpha$ above. We have
\[
       \alpha \in  H^0(\Lambda^{p+1} L_s \otimes \sO(1)) 
                   \cong \Lambda^{p+1} L_s \otimes V
                   \cong \Hom(\Lambda^{p+1} L_s^*,V)
                   \cong_{\sigma'} \Hom(\Lambda^{r-p-1} L,V)
\]
where the last isomorphism is obtained by choosing an orientation 
$\sigma' \in \Lambda^r L_s$. Together with the inclusion
\[
             \iota \colon L_s \hookrightarrow V
\]
we can define
\[
            \pi^* = \iota \oplus \alpha \in \Hom(L_s\oplus \Lambda^{r-p-1} L_s,V).
\]
We denote the induced map on quadrics by the same letter
\[
            \pi^* \in \Hom(L_s\otimes \Lambda^{r-p-1} L_s,S_2(V)).
\]
With this we have a commutative diagram
\begin{center}
\mbox{
\xymatrix{
            {\begin{array}{c} 
                { \scriptstyle \alpha \in } \\
                \Lambda^{p+1} L_s \otimes V
             \end{array}}  
            \ar@<-1ex>[r]
            &
            {\begin{array}{c} 
                { \scriptstyle \sigma \in } \\
                \Lambda^{p} L_s \otimes (I_X)_2
             \end{array}}         
            \\
            {\begin{array}{c} 
                { \scriptstyle \sigma'' \in } \\
                \Lambda^{p+1} L_s \otimes \Lambda^{r-p-1}
             \end{array}}
            \ar[u]^{\id \otimes \alpha}
            \ar@<-1ex>[r]   
            &
            {\begin{array}{c} 
                { \scriptstyle \sigma_{gen} \in } \\
                \Lambda^{p} L_s \otimes (L_s \otimes  \Lambda^{r-p-1})
             \end{array}}
            \ar[u]_{\id \otimes \pi}   
            \\      
            {\begin{array}{c} 
                { \scriptstyle \sigma' \in } \\
                \Lambda^{r} L_s 
             \end{array}} 
            \ar[u]
            \ar[ur]              
         }
    }
\end{center}  
 
where $\sigma'$ maps via $\sigma''$ and $\alpha$ to $\sigma$ by the
construction of $\alpha$. Mapping $\sigma'$ the other way yields a
generic $p$-th syzygy $\sigma_{gen}$ of $\Gensyz_p(L_s)$ by 
proposition/definition \ref{genericsyzygy}. The commutativity of the diagram
shows
\[
             (id \otimes \pi^*)(\sigma_{gen}) = \sigma
\]
Since $\sigma=s$ via the natural isomorphism described above, this 
proves the theorem.
\end{proof}

\begin{cor} \label{notcontained}
With the notations above, $\PZ^n \not\subset \Gensyz_p(L)$.
\end{cor}

\begin{proof} 
Let $s\not=0$ be a syzygy, and
\[
   \pi \colon \PZ(V) \dasharrow \PZ^n
\]
the corresponding linear projection.

Suppose $\PZ^n$ is contained in $\Gensyz_p(L)$. Then by the above
theorem we have
\[
        \Syz(s) = \pi^{-1} (\PZ^n \cap \Gensyz_p(L)) 
                = \pi^{-1} \PZ^n
                = \PZ(V)
\]
Consequently all quadrics involved in $s$ must vanish on $\PZ(V)$,
which is not possible for $s\not=0$.
\end{proof}

\section{$p$-th Syzygies of Rank $p+1$}    \label{reducible}
\nosubsections

The lowest possible rank of a $p$-th linear syzygy is $p+1$. 
As it will turn out, only reducible varieties can have
such a syzygy. To prove this we need:

\begin{prop}
Let $L$ be a $(p+1)$-dimensional vector space. Then
\[
            \Gensyz_p(L) \cong \PZ^p \cup \PZ^0 \subset \PZ^{p+1}
\]
\end{prop}

\begin{proof}
We have $r-p-1=0$, so 
\[
        \Gensyz_p(L) \subset \PZ(L \oplus \Lambda^0 L) \cong \PZ^{p+1}
\]
Let $\{l_1,\dots,l_{p+1}\}$ be a basis of $L$ and $a_0$ a generator
of $\Lambda^0 L$. The ideal of $\Gensyz_p(L)$ is generated by the
the image of
\[
\begin{matrix}
     \Lambda^1 L&\to &L \otimes \Lambda^0 L &\to &S_2(L \oplus \Lambda^0 L) \\
             l_i &\mapsto & l_i \otimes a_0  &\mapsto & l_i\cdot a_0 
\end{matrix}
\]
and consequently
\[
         \Gensyz_p(L) = V(l_1a_0,\dots,l_{p+1}a_0) 
                      = V(l_1,\cdots,l_{p+1}) \cup V(a_0)
                      = \PZ^0 \cup \PZ^p \subset \PZ^{p+1}
\]
\end{proof}

\begin{cor} \label{noreduciblesyzygies}
Let $X \subset \PZ(V)$ be a non degenerate scheme, $I_X$ generated
by quadrics and $s\in V_p$ a $p$-th syzygy of rank $p+1$. Then
$X$ is reducible
\end{cor}

\begin{proof}
Let $L=L_s$ be the space of linear forms involved in $s$.
By theorem \ref{maptogensyz} there exists a linear projection
\[
    \pi \colon \PZ(V) \dasharrow \PZ^n \subset \PZ(L \oplus \Lambda^0 L)
\]
such that
\[
     \pi(X) 
     \subset \pi(\Syz(s)) 
     \subset \PZ^n \cap \Gensyz_p(L)
     = \PZ^n \cap (\PZ^p \cup \PZ^0)
\]
Since $\PZ^n \not\subset \PZ^p$ by corollary \ref{notcontained} and $\PZ^p$ is
a hypersurface in $\PZ^{p+1}$ we have
\[
        \PZ^n \cap \PZ^p = \PZ^{n-1}
\]
Now $X$ is non degenerate in $\PZ(V)$ so $\pi(X)$ is non degenerate in
$\PZ^n$.  Therefore
\[
       \pi(X) \not\subset \PZ^n \cap \PZ^p = \PZ^{n-1}
       \quad \text{and} \quad
       \pi(X) \not\subset \PZ^n \cap \PZ^0. 
\]
Consequently $X$ has to be reducible.
\end{proof}

\begin{defn}
$p$-th syzygies of rank $p+1$ are called reducible syzygies.
\end{defn}

\section{$p$-th Syzygies of Rank $p+2$}      \label{geometric} 
\nosubsections

If $X$ is a non degenerate irreducible variety, the lowest possible
rank of a $p$-th syzygy is $p+2$. As noted by Green and Lazarsfeld
these syzygies are closely connected to linear systems on $X$. We will
recall the corresponding facts in this section.

Lets start by calculating the relevant generic syzygy variety: 
 
\begin{prop}
Let $L$ be a $(p+2)$-dimensional vector space. Then
\[
            \Gensyz_p(L) \cong \PZ^1 \times \PZ^{p+1} \subset \PZ^{2p+3}
\]
where the inclusion is the Segre-embedding.
\end{prop}

\begin{proof}
$r-p-1=1$. Therefore
\[
           \Gensyz_p(L) \subset \PZ(L \oplus \Lambda^1 L) \cong \PZ^{2p+3}
\]
Let $\{l_i\}$ be a basis of $L$ and $\{a_i\}$ a basis of $\Lambda^1 L$.
By proposition \ref{GenSyzEquations} $\Gensyz_p(L)$ is cut out by the image of
\[
\begin{matrix}
      \Lambda^2 L&\to &L \otimes \Lambda^1 L &\to &S_2(L \oplus \Lambda^1 L) \\
       l_i \wedge l_j 
       &\mapsto & l_i \otimes a_j -l_j \otimes a_i  
       &\mapsto & l_ia_j - l_j a_i .
\end{matrix}
\]
Notice that this image is also generated by the $2 \times 2$-minors
of
\[
     M =
     \begin{pmatrix}
     l_1 & \dots & l_{p+2} \\
     a_1 & \dots & a_{p+2} 
     \end{pmatrix}.
\]
Now these are the equations of the Segre embedded $\PZ^1 \times \PZ^{p+1}$.
This proves the proposition.
\end{proof} 

\begin{cor} \label{scrolls}
Let $X \subset \PZ(V)$ be a non degenerate irreducible variety, 
$I_X$ generated by quadrics and $s\in V_p$ a $p$-th syzygy of rank 
$p+2$. Then $X$ is contained in a scroll $\SZ$ of degree $p+2$ and
codimension $p+1$.
\end{cor}

\begin{proof}
Let $L=L_s$ be the space of linear forms involved in $s$.
By theorem \ref{maptogensyz} there exists a linear projection
\[
    \pi \colon \PZ(V) \dasharrow \PZ^n \subset \PZ(L \oplus \Lambda^1 L)
\]
such that
\[
     \pi(X) 
     \subset \pi(\Syz(s)) 
     \subset \PZ^n \cap \Gensyz_p(L)
     = \PZ^n \cap (\PZ^1 \times \PZ^{p+1})
\]
Since $\PZ^1 \times \PZ^{p+1}$ has codimension $p+1$ and degree
$p+2$ in $\PZ(L \oplus \Lambda^1 L)$ we only have to prove that
this intersection is of expected codimension. By Eisenbud 
\cite[Ex. A2.19]{Ei95}
this is the case if the matrix
\[
      M = \begin{pmatrix}
            l_1 & \dots & l_{p+2} \\
            a_1 & \dots & a_{p+2}
          \end{pmatrix}
\]
whose $2 \times 2$-minors cut out $\PZ^1 \times \PZ^{p+1}$ remains
$1$-generic when restricted to $\PZ^n$. Now if $M|_{\PZ^n}$ wasn't
$1$-generic, we would have, after some row and column operations,
a $2 \times 2$-minor of the form
\[
      \det \begin{pmatrix}
           l & l' \\
           0 & a
           \end{pmatrix}
          = l\cdot a
\]
The pullback of this quadric to $\PZ(V)$ is involved in the
syzygy $s$ and therefore contained in the ideal of $X$. But this is
impossible, since $X$ is irreducible and non degenerate.
\end{proof}

This suggests the following definition:

\begin{defn}
The $p$-th syzygies of rank $p+2$ are called scrollar syzygies. The
total space of these syzygies is called the $p$-th space of  
scrollar syzygies.
\end{defn}

We can construct special linear syzygies on $X$ from scrollar syzygies by
intersecting the fibers of the corresponding scroll $\SZ$ with $X$. 
For a canonical curve we have the following well known fact:

\begin{prop}\label{clifford}
Let $C \subset \PZ^{g-1}$ be a non hyperelliptic canonical curve
of genus $g$, $s \in V_p$ a $p$-th scrollar syzygy and 
$D = C \cap \PZ^{g-p-3}$ with $\PZ^{g-p-3}=V(L_s)$.
Then $|D|$ is a special linear system with Clifford index
$\cliff(D) \le g-p-3$.
\end{prop}

\begin{proof}see for example \cite{HCandRanestad}
\end{proof}

\begin{rem}
In \cite{GL1} Green and Lazarsfeld use linear systems $|D|$ of Clifford
index $\cliff(D) = g-p-3$ to construct geometric $p$-th syzygies
of $C$. If $|D|$ is a $g^1_{g-p-1}$ these syzygies are scrollar.
\end{rem}

We can now make a precise statement of the geometric syzygy conjecture
for general canonical curves.

\begin{conj}[Generic Geometric Syzygy Conjecture]
Let $C \subset \PZ^{g-1}$ be a general canonical curve of genus
$g$. Then all minimal rank syzygies are scrollar, and all spaces
of scrollar syzygies are non degenerate.
\end{conj}

\begin{rem}
For special canonical curves it is important to consider the scheme
structure on the space of scrollar syzygies as can be seen in the
case of a curve of genus $6$ with only one $g^1_4$ \cite[p. 174]{AH}.

Also there are geometric $p$-th-syzygies in the sense of Green and Lazarsfeld
\cite{GL1} which are not scrollar. These must also be considered in the
case of special curves. The easiest example of this phenomenon is exhibited 
by the plane quintic curve of genus $6$ \cite{HC}.
\end{rem}

\section{$p$-th Syzygies of Rank $p+3$}    \label{grassmannian}     
\nosubsections

We now consider syzygies whose rank is slightly larger than the rank of
scrollar syzygies. In the remainder of this paper we will see, that
these syzygies imply just enough structure on the minimal free resolution
of a general $K3$ surface with
even sectional genus, to prove the geometric syzygy conjecture for 
canonical curves on these surfaces.

We start again by calculating the relevant generic syzygy variety:  

\begin{prop} \label{grassmannianGenSyz}
Let $L$ be a $(p+3)$-dimensional vector space. Then
\[
            \Gensyz_p(L) \cong \GZ \cup \PZ^{N-p-3} \subset \PZ^N
\]
where $N = { p+4 \choose 2}-1$, $\GZ$ is the Grassmannian $\Gr(p+4,2)$
and the inclusion is the Pl\"ucker embedding.
\end{prop}

\begin{proof}
$r-p-1 = 2$. Therefore 
\[
           \Gensyz_p(L) \subset \PZ(L \otimes \Lambda^2 L) 
                                \cong \PZ^{p+3+{p+3 \choose 2}-1}
                                \cong \PZ^N
\]
Let $\{l_i\}$ be a basis of $L$ and $\{a_{ij}\}$ a basis of $\Lambda^2 L$.
By proposition \ref{GenSyzEquations} $\Gensyz_p(L)$ is cut out by the image of
\[
\begin{matrix}
      \Lambda^3 L&\hookrightarrow 
                 &L \otimes \Lambda^2 L 
                 &\hookrightarrow
                 &S_2(L \oplus \Lambda^2 L) \\
       l_i \wedge l_j \wedge l_k
       &\mapsto & l_i \otimes a_{jk} -l_j \otimes a_{ik} +l_k \otimes a_{ij}
       &\mapsto & l_i a_{jk} -l_j a_{ik} +l_k a_{ij}
\end{matrix}
\]
Notice that these quadrics are also generated by the $4\times4$-Pfaffians
of the skew-symmetric matrix
\[
    M= \begin{pmatrix}
       0 & l_1 &\dots & l_{p+3} \\
    -l_1 & \\
     \vdots & & (a_{ij})\\
    -l_{p+3} 
       \end{pmatrix}
\]
that involve the first row and column. $\GZ$ is cut out by all
$4\times4$-Pfaffians of $M$. Therefore
\[
       \GZ \subset \Gensyz_p(L)
\]
On the other hand the above Pfaffians also vanish, if 
$l_1 = \dots = l_{p+3} =0$. This shows 
\[
       \PZ^{N-p-3} \subset \Gensyz_p(L)
\]
where $\PZ^{N-p-3} = V(l_1, \dots , l_{p+3})$. To prove the proposition
we have to show, that every point of $\Gensyz_p(L)$ outside of
$\PZ^{N-p-3}$ lies on $\GZ$.

Let $x$ be such a point. Since $x \not\in \PZ^{N-p-3}$, there is at
least one linear form of $L$ that doesn't vanish in $x$. Without
restriction we can assume $l_1(x) \not=0$. After the appropriate
row and column transformations $M(x)$ has the form
\[
M(x) =
\begin{pmatrix}
  0      & l_1(x) & 0      & \dots  & 0  \\
 -l_1(x) & 0      & 0      & \dots  & 0 \\
  0      & 0      &        &        &   \\
\vdots   & \vdots &        & M'     &  \\
  0      & 0      &        &        &  \\
\end{pmatrix}
\]
with $M'=(m_{ij})_{i,j\ge3}$ skew symmetric. Now consider the $P$ Pfaffian
that involves the rows and columns $1,2,i$ and $j$:
\[
    P(x) = P_{12ij}(x) = l_1(x)m_{ij}
\]
Since $P$ involves the first row and $x$ is in the generic syzygy scheme,
we have
\[
    l_1(x)m_{ij} = P(x) = 0.
\]
Because $l_1(x)\not=0$, this implies $m_{ij}(x)=0$. Consequently
$M'=0$ and $M(x)$ is of rank $2$. Therefore $x$ is in $\GZ$.
\end{proof}

\begin{cor}
Let $X \subset \PZ(V)$ be a non degenerate irreducible variety, 
$I_X$ generated by quadrics and $s\in V_p$ a $p$-th syzygy of rank 
$p+3$. Then $\pi(X)$ is contained in the Grassmannian $\GZ=\Gr(p+4,2)$.
\end{cor}

\begin{proof}
Let $L=L_s$ be the space of linear forms involved in $s$.
By theorem \ref{maptogensyz} there exists a linear projection
\[
    \pi \colon \PZ(V) \dasharrow \PZ^n \subset \PZ(L \oplus \Lambda^2L)
\]
such that
\begin{align*}
     \pi(X) 
     \subset \pi(\Syz(s)) 
     &= \PZ^n \cap \Gensyz_p(L) \\
     &= \PZ^n \cap (\GZ \cup \PZ^{N-p-3}) \\
     &= (\PZ^n \cap \GZ) \cup (\PZ^n \cap \PZ^{N-p-3}).
\end{align*}
With $\PZ^n \cap \PZ^{N-p-3} \not= \PZ^n$ since $\PZ^n \not\subset
\Gensyz_p(L)$ by corollary \ref{notcontained}.

Now $X$ is non degenerate and irreducible in $\PZ(V)$ and $\pi(X)$ has therefore the same properties in $\PZ^n$. Consequently $\pi(X) \subset \GZ$.
\end{proof}

\begin{defn}
The $p$-th syzygies of rank $p+3$ are called Grassmannian
syzygies.
\end{defn}

\section{Syzygies of $\GZ$}
\nosubsections

We will now study the minimal free resolution of the Grassmannian $\GZ$ 
that occurs in the generic
syzygy variety of a Grassmannian $p$-th syzygy.

I.e. let $\GZ = \Gr(U,2) \subset \PZ(\Lambda^2 U)$ be the 
Grassmannian of $2$-dimensional quotient spaces of the
vector space $U$ with basis $\{u_1,\dots,u_{p+4}\}$.

\begin{prop}
The equations of $\GZ \subset \PZ(\Lambda^2 U)$ are generated by 
the $4 \times 4$-Pfaffians
of
\[
        M_U = \begin{pmatrix}
                0         & u_{12}     & \dots  & u_{1,p+4} \\
               -u_{12}    &     0      & \dots  & u_{2,p+4} \\
               \hdots     & \hdots     & \ddots & \hdots    \\
               -u_{1,p+4} & -u_{2,p+4} & \dots  & 0 \\
              \end{pmatrix}
\]
where the $u_{ij} = u_i \wedge u_j$ are linear forms on $\Lambda^2 U^*$.
\end{prop}

\begin{proof}
For example \cite[Ex. 9.20]{HaAG}.
\end{proof}

\begin{prop} \label{resolutionGrass}
The linear strand of the minimal resolution of $\GZ$ is
\[
      I_{\GZ} \from \Lambda_4 U \otimes \sO(-2)
                   \from \Lambda_{51} U \otimes \sO(-3)
                   \from \dots
                   \from \Lambda_{p+4,1^{p}} U \otimes \sO(-p-2).
\]
\end{prop}

\begin{proof}
The minimal free resolution of an ideal generated by the 
$(2q+2)\times (2q+2)$-Pfaffians of a generic skew symmetric matrix
$M$ is calculated by J\'ozefiak, Pragacz and Weyman in
\cite[Thm 3.14]{PfaffianResolution}.

In our case $q=1$ and we are only interested in the linear strand
of the resolution ($k=1$ in the notation of J\'ozefiak, Pragacz and Weyman).
The $i$th step of the linear strand is then the  Schur functor
corresponding to a Young diagram of the form
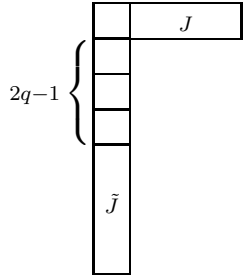
\[
\setlength{\arraycolsep}{-0pt}
\begin{array}{rc}
   \ybox{} & \youngbox{15mm}{\quad {\scriptstyle J} \quad} \\[-1pt]
{\scriptstyle 2q-1}\left\{ \rule{0mm}{4.5ex} \! \right. 
  \begin{matrix} 
    \ybox{} \\[-1pt]
    \ybox{} \\[-1pt] 
    \ybox{} 
  \end{matrix} \\[-1pt]
    \ybox{\rule[-7mm]{0mm}{15mm} \scriptstyle \tilde{J}} \\[-1pt]
\end{array}
\]
where the total number of squares is equal to $2(i+1)$. This proves
the proposition.
\end{proof}

We will now focus our attention on the space $U_p$ of $p$-th syzygies 
of $\GZ$. First we describe the 
syzygies of minimal rank and their spaces of linear forms:

\begin{prop}
The scheme of minimal rank 
$p$-th syzygies of $\GZ$ contains 
the $p$-uple embedding of 
$\PZ^{p+3} \cong \PZ(U^*) =: Y_{min} \hookrightarrow \PZ(U_p^*)$. 
The space of linear forms $L_u$ involved in a minimal rank syzygy 
$u \in U$ is given by 
\[
      L_u = u \wedge U \subset \Lambda^2 U.
\]
The vector bundle of linear forms $\sL$ on $Y_{min}$ is $\sT_{\PZ^{p+3}}(-2)$.
\end{prop}

\begin{proof}
Consider the action of $\GL(p+4)$ of the space of $p$-th syzygies
$U_{p} \cong \Lambda_{p+4,1^p} U$. Now the rank of a syzygy is
invariant under this action and and the space of minimal rank
syzygies is compact. Therefore it has to contain the minimal
orbit
\[
         Y_{min} = G/P \cong \Flag(\PZ^0 \subset \PZ^{p+3}) 
                       \cong \PZ^{p+3} 
                       \xrightarrow{\text{$p$-uple}} 
                       \PZ(U_{p}^*).
\]
Let $u$ be a syzygy in $Y_{min}$.
Since $Y_{min}$ is the minimal orbit, we can without restriction 
assume $u$ to be the maximal weight vector
\[
\setlength{\arraycolsep}{-0pt}    
u =
\begin{array} {cc}
\youngbox{7mm}{1} & \overbrace{\ybox{1} \xyoungbox{\cdots} \ybox{1}}^{p}  \\[-1pt]
\youngbox{7mm}{2} \\[-1pt]
\youngbox{7mm}{\vdots} \\[-1pt]
\youngbox{7mm}{\scriptstyle p+3} \\[-1pt]
\youngbox{7mm}{\scriptstyle p+4} \\[-1pt]
\end{array}
\]
To determine the linear forms involved in $u$ we restrict the
map 
\[
   \psi \colon \Lambda_{p+3,1^{p-1}} U^* 
               \otimes \sO_{\PZ(\Lambda_{p+4,1^{p}} U^*)}(-1)
               \to \Lambda^2 U \otimes 
               \sO_{\PZ(\Lambda_{p+4,1^{p}} U^*)}
\]
from definition \ref{mapofvectorbundles} to the syzygy above. This gives
\[
             \psi|_u = \tildephi(u) \in
                        \Hom(\Lambda_{p+4,1^{p-1}} U^*,\Lambda^2 U)
                        \cong \Lambda_{p+4,1^{p-1}} U \otimes \Lambda^2 U
\]
where
\[
             \tildephi \colon \Lambda_{p+4,1^{p}} U \hookrightarrow
                              \Lambda_{p+4,1^{p-1}} U \otimes \Lambda^2 U
\]
is the map induced by the last step in the linear strand of the
minimal free resolution. Using Young diagrams as in \cite{young} we get:
\[
\setlength{\arraycolsep}{-0pt}    
\begin{array}{ccl}
   \tildephi \colon  & 
\begin{array} {cc}
\youngbox{7mm}{} & \overbrace{\ybox{} \xyoungbox{\cdots} \ybox{} \ybox{}}^{p}  \\[-1pt]
\youngbox{7mm}{} \\[-1pt]
\youngbox{7mm}{\vdots} \\[-1pt]
\youngbox{7mm}{} \\[-1pt]
\youngbox{7mm}{} \\[-1pt]
\end{array}
&
\to
\begin{array} {cc}
\youngbox{7mm}{} & \overbrace{\ybox{} \xyoungbox{\cdots} \ybox{}}^{p-1}  \\[-1pt]
\youngbox{7mm}{} \\[-1pt]
\youngbox{7mm}{\vdots} \\[-1pt]
\youngbox{7mm}{} \\[-1pt]
\end{array}
\otimes
\begin{array} {cc}
\youngbox{7mm}{} \\[-1pt]
\youngbox{7mm}{} \\[-1pt]
\end{array}
\\[60pt]
&
\begin{array} {cc}
\youngbox{7mm}{1} & \ybox{1} \xyoungbox{\cdots} \ybox{1} \ybox{1} \\[-1pt]
\youngbox{7mm}{2} \\[-1pt]
\youngbox{7mm}{\vdots} \\[-1pt]
\youngbox{7mm}{\scriptstyle p+3} \\[-1pt]
\youngbox{7mm}{\scriptstyle p+4} \\[-1pt]
\end{array}
&
\mapsto
\begin{array} {cc}
\youngbox{7mm}{1} & \ybox{1} \xyoungbox{\cdots} \ybox{1}  \\[-1pt]
\youngbox{7mm}{2} \\[-1pt]
\youngbox{7mm}{\vdots} \\[-1pt]
\youngbox{7mm}{\scriptstyle p+3} \\[-1pt]
\end{array}
\otimes
\begin{array} {cc}
\youngbox{7mm}{1} \\[-1pt]
\youngbox{7mm}{\scriptstyle p+4} \\[-1pt]
\end{array}
\pm \dots \pm
\begin{array} {cc}
\youngbox{7mm}{1} & \ybox{1} \xyoungbox{\cdots} \ybox{1}  \\[-1pt]
\youngbox{7mm}{3} \\[-1pt]
\youngbox{7mm}{\vdots} \\[-1pt]
\youngbox{7mm}{\scriptstyle p+4} \\[-1pt]
\end{array}
\otimes
\begin{array} {cc}
\youngbox{7mm}{1} \\[-1pt]
\youngbox{7mm}{2} \\[-1pt]
\end{array}
\end{array}
\]
Consequently the linear forms involved in $u$ are
\[
       L_u = \Img \tildephi(u) 
           = \langle \,
\setlength{\arraycolsep}{-0pt}    
\begin{array} {cc}
\youngbox{7mm}{1} \\[-1pt]
\youngbox{7mm}{\scriptstyle p+4} \\[-1pt]
\end{array}
\, , \dots , \,
\begin{array} {cc}
\youngbox{7mm}{1} \\[-1pt]
\youngbox{7mm}{2} \\[-1pt]
\end{array} \,
             \rangle 
           = \langle u_1\wedge u_{p+4},\dots,u_1 \wedge u_2 \rangle 
           = u \wedge U
\]
as claimed.

We will now determine the vector bundle of linear forms $\sL$ on $Y_{min}$.
For this notice, that the above description of the fibers $L_u$ of $\sL$ 
exhibits $\sL$ as the image of the map
\[
        U \otimes \sO_{\PZ^{p+3}}(-1) \xrightarrow{\wedge U} 
        \Lambda^2 U \otimes \sO_{\PZ^{p+3}}.
\]
which is part of a Koszul complex and factors over $\sT_{\PZ^{p+1}}(-2)$:

\begin{center}
\mbox{
\xymatrix{
             \sO_{\PZ^{p+3}}(-2)
             \ar[r]^-{\wedge U}
             & U \otimes \sO_{\PZ^{p+3}}(-1)
             \ar[rr]^-{\wedge U}
             \ar[dr]
             &&\Lambda^2 U \otimes \sO_{\PZ^{p+3}}
             \ar[r]^-{\wedge U}
             &\Lambda^3 U \otimes \sO_{\PZ^{p+3}}(1) \\
             && \sT_{\PZ^{p+3}}(-2)
                \ar[ur] 
         }
      }
\end{center}

This completes the proof of the proposition.
\end{proof}

\begin{cor}
$\GZ$ has no scrollar $p$-th syzygies.
\end{cor}

\begin{proof}
The proposition above shows $p$-th syzygies have rank greater or equal
to $\rank \sT_{\PZ^{p+3}}= p+3$. Scrollar $p$-th syzygies would
have rank $p+2$.
\end{proof} 

Next we will determine the syzygy varieties of the minimal rank
syzygies:

\begin{prop} \label{minimalsyzvar}
Let $s \in Y_{min}$ be a minimal rank syzygy. Then
\[
       \Syz(s) = \GZ \cup \PZ^{N-p-3} \subset \PZ^N
\]
where $\PZ^{N-p-3} = V(L_s)$ is cut out by the linear forms involved
in $s$.
\end{prop}

\renewcommand{\blackbox}{\ybox{\bullet}}

\begin{proof}
The space $U_{p}$ of $p$-th syzygies of $\GZ$ is isomorphic to 
$\Lambda_{p+4,1^{p}} U$ by proposition \ref{resolutionGrass}. Furthermore
lemma \ref{koszul} exhibits $U_{p}$ as a subspace of
$\Lambda^{p}(\Lambda^2 U) \otimes (I_\GZ)_2$:
\[
\setlength{\arraycolsep}{-0pt}
\begin{array}{rc}
{\scriptstyle p+4}\left\{ \rule{0mm}{7ex} \! \right. 
  \begin{matrix} 
    \ybox{} \\[-1pt]
    \ybox{} \\[-1pt] 
    \ybox{\vdots} \\[-1pt]
    \ybox{} \\[-1pt]
  \end{matrix}
&
\begin{matrix} 
\overbrace{ \xyoungbox{\cdots} \ybox{}}^{p}
\rule[-69pt]{1mm}{0ex} \\[-1pt]
  \end{matrix}
\end{array}
\subset
\left(
\Lambda^{p} \,\,
\begin{array}{c}
    \ybox{} \\[-1pt]
    \ybox{} \\[-1pt] 
\end{array}
\,\right)
\otimes
\begin{array}{c}
    \ybox{} \\[-1pt]
    \ybox{} \\[-1pt] 
    \ybox{} \\[-1pt]
    \ybox{} \\[-1pt] 
\end{array}
= 
\left( 
\Lambda_{\lambda_1} U \oplus \dots \oplus \Lambda_{\lambda_n} U 
\right)
\otimes
\begin{array}{c}
    \ybox{} \\[-1pt]
    \ybox{} \\[-1pt] 
    \ybox{} \\[-1pt]
    \ybox{} \\[-1pt] 
\end{array}
\]
where the $\Lambda_{\lambda_i}$ are the irreducible representations
of $\Lambda^{p+4}(\Lambda^2 U)$. We will now show, that $U_{p}$ is contained
in only one $\Lambda_{\lambda_i} \otimes \Lambda_4$.

To do this, observe, that there are only two ways of coloring $4$
squares of $\Lambda_{p+4,1^{p}}$ which are compatible
with the Littlewood-Richardson 
rule for $\Lambda_{\lambda_i} \otimes \Lambda_4$:
\[
\setlength{\arraycolsep}{-0pt}
\begin{array}{rc}
{\scriptstyle p}\left\{ \rule{0mm}{7ex} \! \right. 
  \begin{matrix} 
    \ybox{} \\[-1pt]
    \ybox{} \\[-1pt] 
    \ybox{\vdots} \\[-1pt]
    \ybox{} \\[-1pt]
  \end{matrix}
&
\begin{matrix} 
\overbrace{ \xyoungbox{\cdots} \ybox{} \ybox{}}^{p}
\rule[-68pt]{1mm}{0ex} 
\end{matrix} \\[-16pt]
\blackbox \\[-1pt]
\blackbox \\[-1pt]
\blackbox \\[-1pt]
\blackbox \\[-1pt]
\end{array}
\text{\quad and \quad}
\begin{array}{rc}
{\scriptstyle p+1}\left\{ \rule{0mm}{9ex} \! \right. 
  \begin{matrix} 
    \ybox{} \\[-1pt]
    \ybox{} \\[-1pt] 
    \ybox{\vdots} \\[-1pt]
    \ybox{} \\[-1pt]
    \ybox{} \\[-1pt]
  \end{matrix}
&
\begin{matrix} 
\overbrace{ \xyoungbox{\cdots} \ybox{} }^{p-1} \! \blackbox
\rule[-81.5pt]{1mm}{0ex} 
\end{matrix} \\[-16pt]
\blackbox \\[-1pt]
\blackbox \\[-1pt]
\blackbox \\[-1pt]
\end{array}
\]

Since the summands of $\Lambda^{p}(\Lambda^2 U)$ have at most
$p$ squares in each row, only the second possibility can occur, i.e
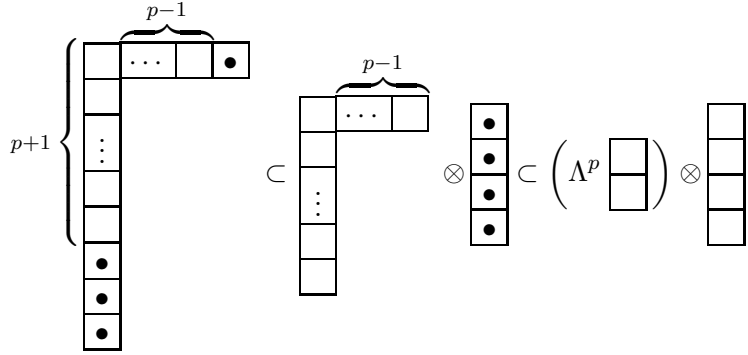
\[
\setlength{\arraycolsep}{-0pt}    
\begin{array}{rc}
{\scriptstyle p+1}\left\{ \rule{0mm}{9ex} \! \right. 
  \begin{matrix} 
    \ybox{} \\[-1pt]
    \ybox{} \\[-1pt] 
    \ybox{\vdots} \\[-1pt]
    \ybox{} \\[-1pt]
    \ybox{} \\[-1pt]
  \end{matrix}
&
\begin{matrix} 
\overbrace{ \xyoungbox{\cdots} \ybox{} }^{p-1} \! \blackbox
\rule[-81.5pt]{1mm}{0ex} 
\end{matrix} \\[-16pt]
\blackbox \\[-1pt]
\blackbox \\[-1pt]
\blackbox \\[-1pt]
\end{array}
\subset \,
\begin{array}{rc}
  \begin{matrix} 
    \ybox{} \\[-1pt]
    \ybox{} \\[-1pt] 
    \ybox{\vdots} \\[-1pt]
    \ybox{} \\[-1pt]
    \ybox{} \\[-1pt]
  \end{matrix}
&
\begin{matrix} 
\overbrace{ \xyoungbox{\cdots} \ybox{} }^{p-1} \! 
\rule[-81.5pt]{1mm}{0ex} 
\end{matrix} \\[-16pt]
\end{array}
\otimes
\begin{matrix}
\blackbox{} \\[-1pt]
\blackbox{} \\[-1pt] 
\blackbox{} \\[-1pt]
\blackbox{} \\[-1pt] 
\end{matrix}
\subset
\left(
\Lambda^{p} \,\,
\begin{array}{c}
    \ybox{} \\[-1pt]
    \ybox{} \\[-1pt] 
\end{array}
\,\right)
\otimes
\begin{array}{c}
    \ybox{} \\[-1pt]
    \ybox{} \\[-1pt] 
    \ybox{} \\[-1pt]
    \ybox{} \\[-1pt] 
\end{array}
\]
By lemma \ref{quadricsinvolved} the syzygy scheme $\Syz(s)$ is cut out by the
quadrics involved in the image of $s$ under the above inclusion.
Without loss of generality we can assume $s$ to be the maximal weight vector
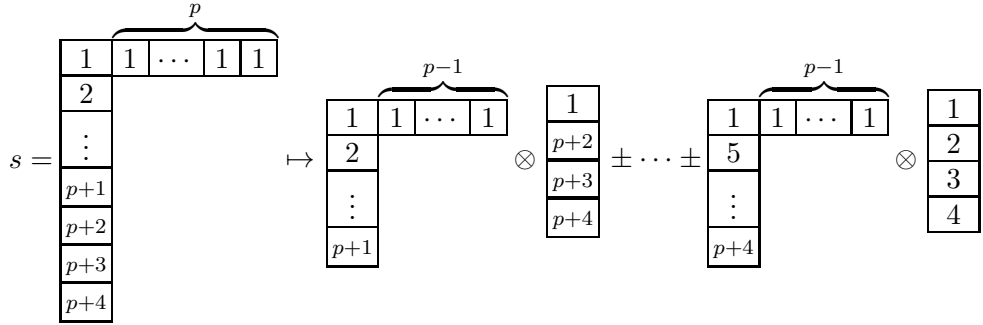
\[
\setlength{\arraycolsep}{-0pt}    
s =
\begin{array} {cc}
\youngbox{7mm}{1} 
& \overbrace{\ybox{1} \xyoungbox{\cdots} \ybox{1} 
             \ybox{1}}^{p}  
\\[-1pt]
\youngbox{7mm}{2} \\[-1pt]
\youngbox{7mm}{\vdots} \\[-1pt]
\youngbox{7mm}{\scriptstyle p+1} \\[-1pt]
\youngbox{7mm}{\scriptstyle p+2} \\[-1pt]
\youngbox{7mm}{\scriptstyle p+3} \\[-1pt]
\youngbox{7mm}{\scriptstyle p+4} \\[-1pt]
\end{array}
\mapsto\,
\begin{array} {cc}
\youngbox{7mm}{1} & \overbrace{\ybox{1} \xyoungbox{\cdots} \ybox{1} }^{p-1}  \\[-1pt]
\youngbox{7mm}{2} \\[-1pt]
\youngbox{7mm}{\vdots} \\[-1pt]
\youngbox{7mm}{\scriptstyle p+1} \\[-1pt]
\end{array}
\otimes\,
\begin{array} {cc}
\youngbox{7mm}{1} \\[-1pt]
\youngbox{7mm}{\scriptstyle p+2} \\[-1pt]
\youngbox{7mm}{\scriptstyle p+3} \\[-1pt]
\youngbox{7mm}{\scriptstyle p+4} \\[-1pt]
\end{array}
\, \pm \dots \pm
\begin{array} {cc}
\youngbox{7mm}{1} & \overbrace{\ybox{1} \xyoungbox{\cdots} \ybox{1} }^{p-1}  \\[-1pt]
\youngbox{7mm}{5} \\[-1pt]
\youngbox{7mm}{\vdots} \\[-1pt]
\youngbox{7mm}{\scriptstyle p+4} \\[-1pt]
\end{array}
\otimes\,
\begin{array} {cc}
\youngbox{7mm}{1} \\[-1pt]
\youngbox{7mm}{2} \\[-1pt]
\youngbox{7mm}{3} \\[-1pt]
\youngbox{7mm}{4} \\[-1pt]
\end{array}
\]
Consequently the syzygy scheme $\Syz(s)$ of $s$ is cut out by
the $4 \times 4$-Pfaffians that involve the first row and
column of $M$. The same argument as in the proof of proposition 
\ref{grassmannianGenSyz} shows
\[
         \GZ \cup \PZ^{N-p-3} = \Syz(s)
\]
\end{proof}

If $s \in U_{p}$ is any $p$-th syzygy, the syzygy scheme
might be more complicated, but its ideal still contains 
certain special quadrics:

\begin{defn}
A quadric $Q$ is called a generalized $4 \times 4$-Pfaffian of a skew symmetric
matrix $M$, if there exists an invertible matrix $B$ such that
$Q$ is a $4\times4$-Pfaffian of $B^tMB$.
\end{defn}

\begin{rem}
In our case the generalized $4 \times 4$-Pfaffians 
of $M_U$ are the decomposable
elements in $(I_\GZ)_2 = \Lambda_4 U$. 
\end{rem}

\begin{lem}\label{generalisedpaffian}
Let $s \in U_{p}$ a $p$-th syzygy of $\GZ$. Then the ideal
of $\Syz(s)$ contains a generalized $4 \times 4$ Pfaffian.
\end{lem}

\begin{proof}
Recall the inclusion
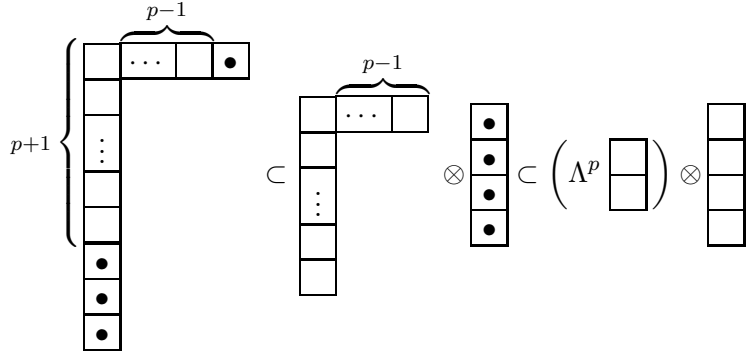
\[
\setlength{\arraycolsep}{-0pt}    
\begin{array}{rc}
{\scriptstyle p+1}\left\{ \rule{0mm}{9ex} \! \right. 
  \begin{matrix} 
    \ybox{} \\[-1pt]
    \ybox{} \\[-1pt] 
    \ybox{\vdots} \\[-1pt]
    \ybox{} \\[-1pt]
    \ybox{} \\[-1pt]
  \end{matrix}
&
\begin{matrix} 
\overbrace{ \xyoungbox{\cdots} \ybox{} }^{p-1} \! \blackbox
\rule[-81.5pt]{1mm}{0ex} 
\end{matrix} \\[-16pt]
\blackbox \\[-1pt]
\blackbox \\[-1pt]
\blackbox \\[-1pt]
\end{array}
\subset \,
\begin{array}{rc}
  \begin{matrix} 
    \ybox{} \\[-1pt]
    \ybox{} \\[-1pt] 
    \ybox{\vdots} \\[-1pt]
    \ybox{} \\[-1pt]
    \ybox{} \\[-1pt]
  \end{matrix}
&
\begin{matrix} 
\overbrace{ \xyoungbox{\cdots} \ybox{} }^{p-1} \! 
\rule[-81.5pt]{1mm}{0ex} 
\end{matrix} \\[-16pt]
\end{array}
\otimes
\begin{matrix}
\blackbox{} \\[-1pt]
\blackbox{} \\[-1pt] 
\blackbox{} \\[-1pt]
\blackbox{} \\[-1pt] 
\end{matrix}
\subset
\left(
\Lambda^{p} \,\,
\begin{array}{c}
    \ybox{} \\[-1pt]
    \ybox{} \\[-1pt] 
\end{array}
\,\right)
\otimes
\begin{array}{c}
    \ybox{} \\[-1pt]
    \ybox{} \\[-1pt] 
    \ybox{} \\[-1pt]
    \ybox{} \\[-1pt] 
\end{array}
\]
from the last proof. For simplicity we will call the first 
Young diagram a big hook and
the second one a small hook. The third one we will call a line.

Now $U_{p}$ has a basis $\{s_\beta\}$ enumerated by big hooks of
the form
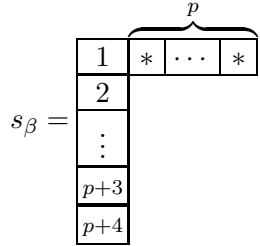
\[
\setlength{\arraycolsep}{-0pt}    
s_\beta =
\begin{array} {cc}
\youngbox{7mm}{1} & \overbrace{\ybox{*} \xyoungbox{\cdots} \ybox{*}}^{p}  \\[-1pt]
\youngbox{7mm}{2} \\[-1pt]
\youngbox{7mm}{\vdots} \\[-1pt]
\youngbox{7mm}{\scriptstyle p+3} \\[-1pt]
\youngbox{7mm}{\scriptstyle p+4} \\[-1pt]
\end{array}
\]
The image of such a big hook $s_\beta$ under the above map is
a sum of tensor products of small hooks and lines, where each
product contains the same entries as $s_\beta$. Notice that the small hook
$d$ of maximal weight can therefore only 
occur in the image of the the following big hooks:
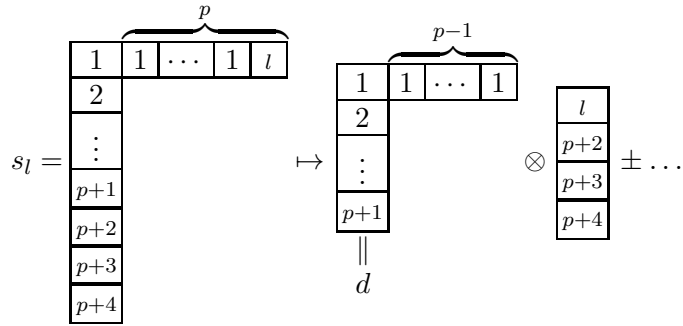
\[
\setlength{\arraycolsep}{-0pt}    
s_l =
\begin{array} {cc}
\youngbox{7mm}{1} 
& \overbrace{\ybox{1} \xyoungbox{\cdots} \ybox{1} 
             \ybox{\scriptstyle l}}^{p}  
\\[-1pt]
\youngbox{7mm}{2} \\[-1pt]
\youngbox{7mm}{\vdots} \\[-1pt]
\youngbox{7mm}{\scriptstyle p+1} \\[-1pt]
\youngbox{7mm}{\scriptstyle p+2} \\[-1pt]
\youngbox{7mm}{\scriptstyle p+3} \\[-1pt]
\youngbox{7mm}{\scriptstyle p+4} \\[-1pt]
\end{array}
\mapsto\,
\begin{array} {cc}
\youngbox{7mm}{1} & \overbrace{\ybox{1} \xyoungbox{\cdots} \ybox{1} }^{p-1}  \\[-1pt]
\youngbox{7mm}{2} \\[-1pt]
\youngbox{7mm}{\vdots} \\[-1pt]
\youngbox{7mm}{\scriptstyle p+1} \\[2pt]
\parallel \\
d\\
\end{array}
\otimes\,
\begin{array} {cc}
\youngbox{7mm}{\scriptstyle l} \\[-1pt]
\youngbox{7mm}{\scriptstyle p+2} \\[-1pt]
\youngbox{7mm}{\scriptstyle p+3} \\[-1pt]
\youngbox{7mm}{\scriptstyle p+4} \\[-1pt]
\end{array}
\, \pm \dots
\]
Let now $s = \sum \mu_\beta s_\beta$ be any $p$-th syzygy of $\GZ$, i.e a 
linear combination of big hooks. 
Then the image of $s$ can be written as $\sum d_\alpha \otimes Q_\alpha$ with
$\{ d_\alpha \}$ a basis of $\Lambda_{p+1,1^{p-1}} U$ enumerated by
small hooks. By the argument above the small hook 
$d$ of maximal weight occurs with a quadric $Q$ of the form
\[
    Q = \sum \mu_{l} \,
\setlength{\arraycolsep}{-0pt}    
\begin{array} {c}
\youngbox{7mm}{\scriptstyle l} \\[-1pt]
\youngbox{7mm}{\scriptstyle p+2} \\[-1pt]
\youngbox{7mm}{\scriptstyle p+3} \\[-1pt]
\youngbox{7mm}{\scriptstyle p+4} \\[-1pt]
\end{array}
\,
= \sum \mu_l (u_l \wedge u_{p+2} \wedge u_{p+3} \wedge u_{p+4})
= \left(  \sum \mu_l u_l \right) \wedge u_{p+2} \wedge u_{p+3} \wedge u_{p+4}
\]
Consequently $Q$ is either zero or a generalized $4 \times 4$-Pfaffian.

We will now check, that we can assume $Q\not\equiv 0$. For this consider
the minimal orbit $G/P$ of $\GL(p+4)$ in the space $\Lambda_{p+1,1^{p-1}} U$
of small hooks. Since this representation is irreducible, $G/P$ is
non degenerate. We can therefore choose a basis $\{ d_j \}$ of this space
consisting only of points in the minimal orbit. Now write the image
of $s$ in this basis:
\[
       s \mapsto \sum_j d_j \otimes Q_j
\]
Since $s \not\equiv 0$, there is at least one $Q_j \not\equiv 0$.
Since the corresponding $d_j$ is in the minimal orbit of $\GL(p)$
we can after a coordinate change of $U$ assume it to be the
maximal weight vector $d$ above. Applying the above argument to
the transformed $Q_j \not\equiv 0$ shows that $Q_j$ is a generalized
$4 \times 4$-Pfaffian.
\end{proof}

\begin{rem}
The argument above even shows, that the ideal of
$\Syz(s)$ is generated by generalized Pfaffians.
\end{rem}

\section{Linear sections of $\GZ$}
\nosubsections

Let $X \subset \PZ(V)$ be a non degenerate irreducible variety,
$s \in V_p$ a Grassmannian syzygy and
\[
       \pi \colon \PZ(V) \dasharrow \PZ^n \subset \PZ(\Lambda^2 U)
\]
the induced projection with
\[
         \pi(X) \subset \PZ^n \cap \GZ.
\]
In this situation we have a natural map from $p$-th syzygies of $\GZ$ to
$p$-th syzygies of $X$
\[
       \alpha_p \colon H^0\bigl(\PZ(\Lambda^2 U),
                                \Omega^p(p+2) \otimes I_{\GZ} \bigr)
                \to    H^0\bigl(\PZ(V),\Omega^p(p+2) \otimes I_{X} \bigr)
\]
given by restriction of the corresponding twisted $p$-forms to $\PZ^n$
and pulling them back to $\PZ(V)$. If a syzygy $s$ is not in the
kernel of $\alpha_p$ we say ``s survives the restriction to $X$''.

In this section we want to prove that $\alpha_{p}$ is injective, i.e.
that all $p$-th syzygies of $\GZ$ survive the restriction to 
$X$.

The first step is

\begin{prop}
All minimal rank syzygies $u \in Y_{min}$ survive the restriction
to $X$.
\end{prop}

\begin{proof}
By theorem \ref{maptogensyz} there exist a generic $p$-th syzygy 
$u_{gen} \in Y_{min}$ of 
\[
        \Gensyz_{p}(L_s) = \GZ \cup \PZ_s
\]
with
\[
             \PZ_s = V(L_s)
\]
and 
\[
            s = \pi^*(u_{gen}).
\]
Consequently $u_{gen}$ survives the restriction to $X$. Furthermore
we have $L_s = L_{u_{gen}}$, i.p. no linear form in $L_{u_{gen}}$
vanishes on $\PZ^n$.

Suppose $u\in Y_{min}$ is a minimal rank syzygy that doesn't survive
the restriction to $X$, i.e.
\[
            \PZ^n \subset \Syz(u)=\GZ \cup \PZ_u
\]
where $\PZ_u = V(L_u)$. Now $\PZ^n \not\subset \GZ$ so we must have 
$\PZ^n \subset \PZ_u$. In particular all linear forms in $L_u$
must vanish on $\PZ^n$. But this is impossible, since
$u \wedge u_{gen}$ is a linear form in $L_u$ and $L_{u_{gen}}$
that doesn't vanish on $\PZ^n$ by the above argument.
\end{proof}

To extend this result to arbitrary $p$-th syzygies, we need

\begin{prop} \label{pfaffianssurvive}
Let $X$ be an irreducible variety as above. Then
all generalized Pfaffians of $M_U$ survive the restriction to $X$
\end{prop}

\begin{proof}
Suppose $P$ is a generalized Pfaffian that doesn't survive
the restriction to $X$, i.e. $\PZ^n \subset V(P)$. Without loss
of generality we can assume
\[
      P = u_1 \wedge u_2 \wedge u_3 \wedge u_4 
        = \pfaff(M_{1234})
        = \pfaff 
          \begin{pmatrix}
                   0 & u_{12}  & u_{13}  & u_{14} \\
             -u_{12} & 0       & u_{23}  & u_{24} \\
             -u_{13} & -u_{23} & 0       & u_{34} \\
             -u_{14} & -u_{24} & -u_{34} & 0      
          \end{pmatrix}
\]
First of all we will prove that $\PZ^n \subset V(P)$ implies
that the restriction of $M_{1234}$ to $\PZ^n$ can, after a 
suitable coordinate change, be written in the form
\[
     M_{1234}|_{\PZ^n}=
         \begin{pmatrix}
                0 & 0 & 0 & * \\
                0 & 0 & * & * \\
                0 & * & 0 & * \\
                * & * & * & 0
         \end{pmatrix}.
\]
To see this consider the vertex $\PZ^{N-6}=V(u_{12},\dots,u_{34})$
of $V(P)$ and project from there:
\[
       \phi \colon \PZ(\Lambda^2 U) \dasharrow \PZ^5
\]   
The image of $V(P)$ is the Grassmannian $\Gr(2,4)$. $\Gr(2,4)$ is
also cut out by the Pfaffian $P$. Now $\PZ^n \subset V(P)$ so either
$\PZ^n \subset \PZ^{N-6}$ and
\[
     M_{1234}|_{\PZ^n}=
         \begin{pmatrix}
                0 & 0 & 0 & 0 \\
                0 & 0 & 0 & 0 \\
                0 & 0 & 0 & 0 \\
                0 & 0 & 0 & 0
         \end{pmatrix}
\]
or $\phi(\PZ^n) \subset \Gr(2,4)$. Since $\Gr(2,4)$ is a quadric of
dimension $4$, $\dim \phi(\PZ^n) \le 2$ \cite[p. 735]{GH}.

If $\phi(\PZ^n) \cong \PZ^0$ is a point in $\Gr(2,4)$, $M_{1234}|_{\PZ^0}$ is a
matrix of rank $2$. Therefore after a suitable coordinate change we 
have
\[
     M_{1234}|_{\PZ^n}=
         \begin{pmatrix}
                0 & 0 & 0 & 0 \\
                0 & 0 & 0 & 0 \\
                0 & 0 & 0 & * \\
                0 & 0 & * & 0
         \end{pmatrix}.
\]
If $\phi(\PZ^n) \cong \PZ^1$ is a line in $\Gr(2,4)$, then this
line is a Schubert-cycle
\[
             \PZ^1 = \{ l \in \Gr(2,4) \,| \, p_0 \in l \subset H_0\}
\]
with $p_0 \in \PZ^3$ an point and $H_0 \subset \PZ^3$ a
hyperplane \cite[p. 757]{GH}. Consequently
\[
     M_{1234}|_{\PZ^n}=
         \begin{pmatrix}
                0 & 0 & 0 & 0 \\
                0 & 0 & 0 & * \\
                0 & 0 & 0 & * \\
                0 & * & * & 0
         \end{pmatrix}
\]
after a suitable coordinate change.

If $\phi(\PZ^n) \cong \PZ^2$ lies in $\Gr(2,4)$, then there
are two possibilities. Firstly
\[
      \PZ^2 = \{ l \in \Gr(2,4) \, | \, l \subset H_0 \}
\]
with $H_0 \subset \PZ^3$ a hyperplane, and
\[
     M_{1234}|_{\PZ^n}=
         \begin{pmatrix}
                0 & 0 & 0 & 0 \\
                0 & 0 & * & * \\
                0 & * & 0 & * \\
                0 & * & * & 0
         \end{pmatrix}
\]
or
\[
      \PZ^2 = \{ l \in \Gr(2,4) \, | \, p_0 \in l\}
\]
with $p_0 \in \PZ^3$ a point, and
\[
     M_{1234}|_{\PZ^n}=
         \begin{pmatrix}
                0 & 0 & 0 & * \\
                0 & 0 & 0 & * \\
                0 & 0 & 0 & * \\
                * & * & * & 0
         \end{pmatrix}
\]
This proves out claim about the shape of $M_{1234}|_{\PZ^n}$. For the
whole matrix $M_U$ we have therefore, after a coordinate change
\[
         M_U|_{\PZ^n} = \begin{pmatrix}
                        0 & 0 & 0 & * & \dots & * \\
                        0 \\
                        0 \\
                        * &   &   & * \\
                        \vdots \\
                        * \\
                        \end{pmatrix}
\]
Now the linear forms in the first row of $M_U$ are the
linear forms of $L_u$ for a particular $u \in Y_{min}$. During the
restriction to $X$ at least two of these linear forms vanish
because of the shape of $M_U|_{\PZ^n}$. Therefore
\[
           \rank u|_X \le p+3-2 = p+1.
\]
But this is impossible for irreducible $X$ by corollary 
\ref{noreduciblesyzygies}.
Consequently $\PZ^n \not\subset V(P)$ for all generalized Pfaffians
$P$ of $M_U$.
\end{proof}

\begin{cor} \label{injective}
If $X$ is irreducible, then the restriction $\alpha_{p}$ 
of $p$-th syzygies of $\GZ$ to
$p$-th syzygies of $X$ is injective.
\end{cor}

\begin{proof}
Let $s\in V_{p}$ be a $p$-th syzygy of $\GZ$. By lemma
\ref{generalisedpaffian} the ideal of $\Syz(s)$ contains at least one
generalized Pfaffian $P$ of $M_U$. Since $P$ survives the
restriction to $X$ by proposition \ref{pfaffianssurvive}, 
$\Syz(s)$ cannot contain
$\PZ^n$. Consequently $s$ also survives the restriction. 
\end{proof}

Finally we will describe which syzygies
$u \in Y_{min}$ of $\GZ$ drop rank during the restriction to $X$. For this let
$\PZorth \subset \PZ(\Lambda^2 U^*)$ be the orthogonal space of $\PZ^n$
and $\GZ^* = \Gr(U^*,2)$ the dual Grassmannian.

\begin{lem} \label{droprank}
Let
\[
     u \in \GZ^* \cap \PZorth 
\]
be a decomposable linear form in $\PZorth$.

Then all syzygies in the line of $Y_{min} \cong \PZ(U^*)$ 
corresponding to $u=u' \wedge u''$ 
drop rank when restricted to $X$.
\end{lem}

\begin{proof}
Consider a syzygy $s \in Y_{min}$, and its space of
linear forms $L_s \subset \PZ(\Lambda^2 U^*)$. 

When we restrict $s$ to $\PZ^n$ all 
linear forms $l$ in $\PZorth \cap L_s$ vanish. So all syzygies whose
space of linear forms $L_s$ intersects $\PZorth$ drop rank when restricted
to $X$. 

In our case consider the line spanned by the
minimal rank syzygies $u'$ and $u''$ in $Y_{min} \cong \PZ(U^*)$.
The space of linear forms of a syzygy $\lambda u' + \mu u''$
on this line is  
\[
    L_{\lambda:\mu} = (\lambda u' + \mu u'') \wedge U
\]
and therefore contains $u = u' \wedge u'' \subset \PZorth$. 
Consequently all syzygies on this
line drop rank during restriction.
\end{proof}

\section{Syzygies of $S$}
\nosubsections

Let $S \subset \PZ^{g}$ be a $K3$-Surface whose Picard group is generated by
$\sO(C)$ where $C$ is a smooth curve of even genus $g=2k$.

We first prove some standard facts about $S$:

\begin{lem}
The ideal of $S$ contains no rank $4$ quadric.
\end{lem}

\begin{proof}
Suppose $Q$ is a rank $4$ quadric with $S \subset Q$. Then the rulings
of $Q$ cut out two linear series $|C_1|$ and $|C_2|$ with 
$|C_1 + C_2|=|C|$. This is impossible, since the Picard group
of $S$ is generated by $|C|$.
\end{proof}

\begin{cor} \label{nogeosyzforK3}
$S$ has no scrollar syzygies.
\end{cor}

\begin{proof}
By corollary \ref{scrolls} a scrollar syzygy implies the existence
of a scroll $\SZ$ containing $S$. These scrolls are cut out
by rank $4$ quadrics. Since $S$ is contained in no rank $4$
quadric this is impossible.
\end{proof} 

\begin{cor}
$S$ has a Grassmannian $(k-2)$-nd syzygy
\end{cor}

\begin{proof}
Consider a general hyperplane section $C \cap \PZ^{g-1}$ with
$\PZ^{g-1} = V(l)$. Since $S$ is arithmetically Cohen Macaulay, 
the restriction maps
\[
       \alpha_p \colon H^0\bigl(\PZ^g,
                                \Omega^p(p+2) \otimes I_{S} \bigr)
                \to    H^0\bigl(\PZ^{g-1},\Omega^p(p+2) \otimes I_{C} \bigr)
\]
are isomorphisms by \cite[Thm 3.b.7]{GreenKoszul}.

Let $s$ be a $p$-th syzygy of $S$ and 
\[
      \tildephi(s) \colon V_{p-1}^* \to V
\]
the map from definition \ref{tildephi}, 
where $V_{p-1}$ is the space of $(p-1)$-st syzygies of $S$, and
$\PZ(V)=\PZ^g$. Then $\Img(\tildephi(s))=L_s$. Restriction to $C$
gives the diagram
\begin{center}
\mbox{
\xymatrix{
           V_{p-1} \ar[rr]^{\tildephi(s)} \ar@{=}[d]_{\alpha_{p-1}}
           && V \ar@{->>}[d]^{\iota^*}\\     
           V_{p-1}  \ar[rr]^{\tildephi(\alpha_p(s))}
           &&
           V'
         }
     }
\end{center}
where $\PZ(V') = \PZ^{g-1}$, and $\iota^*$ is the natural projection
induced by the inclusion 
\[
    \iota \colon \PZ^{g-1} \hookrightarrow \PZ^{g}.
\]
The kernel of $\iota^*$ is generated by $l$. The space of linear forms
involved in $\alpha_p(s)$ is therefore
\[
     L_{\alpha_p(s)} = \Img(\tildephi(\alpha_p(s)) 
                     = \iota^* (\Img(\tildephi(s))
                     = \iota^* L_s
\]
In particular we have
\[
       \rank \alpha_p(s) = \left\{
       \begin{matrix}
          \rank s -1 & \text{if $l \in L_s$} \\
          \rank s    & \text{if $l \not\in L_s$} 
       \end{matrix}
       \right.
\]
Now $C$ has at least finitely many $g^1_{k+1}$'s which induce
scrollar $(k-2)$-nd syzygies by the construction of 
Green and Lazarsfeld \cite{GL1}.

By prop \ref{nogeosyzforK3} these cannot come from scrollar syzygies of
$S$, so they have to come from syzygies with higher rank. Since
the difference in rank can be at most $1$ by the above argument
these syzygies are Grassmannian.
\end{proof}

Consider now such a Grassmannian $(k-2)$-nd syzygy of $S$. It induces
a linear projection
\[
    \pi \colon \PZ^g \dasharrow \PZ^n \subset \PZ^N
\]
with $N = {k+2 \choose 2}-1$ and
\[
    \pi(S) \subset \GZ = \Gr(k+2,2) \subset \PZ^N.
\]
We even have

\begin{prop}
$\pi$ is an embedding, i.e. $\PZ^g \cong \PZ^n$ via $\pi$.
\end{prop}

\begin{proof}
Consider the orthogonal space $\PZorth$ of $\PZ^n$ and the
dual Grassmannian $\GZ^*=\Gr(2,k+2)$. Since $S$ has no scrollar
syzygies, $\PZorth$ does not intersect $\GZ^*$ by proposition \ref{droprank}.

Now $\dim \GZ^*= 2k =g$, and consequently
\[
      1 + \dim \PZ^n = \codim \PZorth \ge \dim \GZ^* + 1 = g+1.
\]
Since $\pi$ is surjective on $\PZ^n$ this proves the proposition.
\end{proof}

\begin{rem}
Notice that this is the embedding constructed by Mukai in 
\cite{ClassificationMukai}
with vector bundle methods. His rank $2$ vector bundle $E$ is the
restriction of the universal quotient bundle $Q$ on $\GZ$ 
to $S\cong \pi(S) \subset \GZ$.

Notice also that our construction also gives an algorithm to determine $\pi$ 
explicitly from a Grassmannian syzygy $s$ by lifting $\tildephi(s)$ to a 
map of complexes as in theorem \ref{maptogensyz}.
\end{rem}

Using the theorem of Voisin we can now describe the 
$(k-2)$-nd syzygies of $S$ geometrically:

\begin{prop}
The space of minimal rank syzygies of $S$
contains a $(k-2)$-uple embedded $\PZ^{k+1}$.
Further more the space of all $(k-2)$-nd syzygies of $S$ 
is isomorphic to the ambient space of this embedding.
\end{prop}

\begin{proof}
Let $U_{k-2}$ be the space of $(k-2)$-nd syzygies of $\GZ = \Gr(U,2)$ and
$V_{k-2}$ the corresponding space of $(k-2)$-nd syzygies of $S$.

By the corollary \ref{injective}  the map
\[
    \alpha_{k-2} \colon U_{k-2} \to V_{k-2}
\]
induced by restriction of syzygies is injective. Consequently
\[
       \dim V_{k-2} \ge \dim U_{k-2} = \dim \Lambda_{k+2,1^{k-2}} U 
                                     = \dim S^{k-2} U =
       { 2k-1 \choose k-2 }.
\]
On the other hand, the Hilbert function
of $S$ gives:
\[
       \dim V_{p} - \beta_{p,p+2} 
       = (p+1) { g-2 \choose p+2} - (g-p+2){g-2 \choose g-p-1}
\]
(see for example \cite{Sch91}). In our case $p=k-2$, $g=2k$ and
$\beta_{p,p+2}=0$ by Voisin's theorem. Consequently
\[
\dim V_{k-2}
       = (k-1) { 2k-2 \choose k} - k{2k-2 \choose k+1}  
       = { 2k-1 \choose k+1 }    
       = { 2k-1 \choose k-2 }    
\]
and $V_{k-2} \cong U_{k-2}$. 

Since $\PZorth$ doesn't intersect $\GZ^*$ all syzygies in 
\[
       Y_{min} \cong \PZ^{k+1} \xrightarrow{\text{$(k-2)$-uple}} \PZ(U_{k-2}^*)
\]
remain of rank $k+1$ during restriction. They also remain  minimal,
since $S$ has no scrollar $(k-2)$-nd syzygies of rank $k$ by 
corollary \ref{nogeosyzforK3}. 
\end{proof}

\section{Syzygies of $C$}
\nosubsections

Now consider a general linear section
\[
        C = S \cap \PZ^{g-1}.
\]
During this further restriction some syzygies drop rank and become
scrollar syzygies:

\begin{prop}
The scrollar $(k-2)$-nd syzygies of $C$ from a configuration of
$\frac{1}{k+1}{ 2k \choose k}$ lines in $\PZ^{k+1}$ that are
embedded in the space of $(k-2)$-nd syzygies of $C$ as 
rational normal curves of degree $(k-2)$ on a $(k-2)$-uple
embedding of $\PZ^{k+1}$.
\end{prop} 

\begin{proof}
Since
$C$ is a general linear section of $S$, their spaces of syzygies are
isomorphic. We now determine which syzygies do drop rank. With lemma
\ref{droprank} we find some of these, when the orthogonal space 
$\PZorth_C$ of $\PZ^{g-1}$ intersects $\GZ^*$.

$\PZorth_C$ contains
the orthogonal space $\PZorth_S$ of $\PZ^{g}$. Since $\PZorth_S$ doesn't
intersect $\GZ^*$, $\PZorth_C$ can only intersect in finitely many points.
On the other hand $\dim \GZ^* = 2k = g = \codim \PZorth_C$ so this is
also the expected intersection dimension. Consequently the
number $r$ of intersection points is equal to the degree of $\GZ^*$.
A formula for the degree of Grassmannians can be found in 
\cite[p. 247]{HaAG}:
\[
         r = \deg \GZ^* = (2k)! \prod_{i=0}^1 \frac{i!}{(k+i)!}
                        = \frac{(2k)!}{k!(k+1)!}
                        = \frac{1}{k+1}{ 2k \choose k}
\]
Now each point of the intersection corresponds to a line
of scrollar syzygies in $\PZ^{k+1} \cong \PZ(U^*)$. This gives
our configuration of lines.

To prove these are all scrollar syzygies of $C$, recall that the
syzygy variety of a scrollar $(k-2)$-nd syzygy is a scroll $\SZ$ whose
fibers cut out linear equivalent divisors. These divisors are part
of a linear system with Clifford index 
\[
        \cliff(D) \ge g-(k-2)-3=k-1
\]
by proposition \ref{clifford}.

Since $C$ is general in the sense of Brill-Noether-Theory by
\cite{BrillNoeterLazarsfeld}, 
the only linear systems on $C$ with Clifford index $k-1$
are the $g^1_{k+1}$'s. In other words: all scrollar $(k-2)$-nd 
syzygies of $C$ come from $g^1_{k+1}$'s. 

Now each $g^1_{k+1}$ induces a $\PZ^1$ of scrollar 
$(k-2)$-nd syzygies. See for example \cite[Lem 2.2.8]{HCandRanestad} 
for a proof of this elementary
fact. 

A formula for the number of $g^1_{k+1}$'s in  $W^1_{k+1}(C)$ 
for a Brill-Noether-general curve $C$
is given in \cite[p. 211]{ACGH}:
\[
         \deg W^1_{k+1} = g! \prod_{i=0}^1 \frac{i!}{(g-(k+1)-1+i)!}
                        = 2k! \prod_{i=0}^1 \frac{i!}{(k+i)!}
                        = \deg \GZ^*
\]
So there are no  scrollar $(k-2)$-nd syzygies of $C$ except 
the ones in the configuration above.
\end{proof}
 
We now deduce the desired special case of the geometric syzygy conjecture

\begin{thm}
Let $C$ be a general hyperplane section of a $K3$ surface $S$ whose
Picard group is generated by $C$. Then the space of $(k-2)$-nd 
scrollar syzygies of $C$ is non degenerate.
\end{thm}

\begin{proof}
Let $Z$ the configuration of lines in $\PZ^{k+1}$ that correspond
to scrollar syzygies by the previous proposition. 
The $(k-2)$-uple embedding of $\PZ^{k+1}$ is
\[
          \iota \colon \PZ^{k+1} \hookrightarrow \PZ^N
\]
with $N={2k-1 \choose k-2}-1$. We have to show, that $\iota(Z)$
spans $\PZ^N$.

For this notice that $Z$ is
the locus of syzygies $s$ where the dimension of the space of linear forms 
$L_s$ is $k$. $Z$ can therefore be described by the determinantal locus where
the composition map of vector bundles
\[
       \beta \colon \sL \to \Lambda^2 U \otimes \sO_{\PZ^{k+1}}
                      \to V' \otimes \sO_{\PZ^{k+1}}
\]
drops rank. Here $V'$ denotes the space of linear forms on $\PZ^{g-1}$ and
$\sL \cong \sT_{\PZ^{k+1}}(-2)$ the vector bundle 
of linear forms on $Y_{min} \cong \PZ^{k+1}$. 
Notice that $\beta$ drops rank
in expected dimension
\[
     \dim \PZ^{k+1} - (\rank \sL - k)(\dim V' - k) = k+1 - (k+1-k)(2k-k) = 1
\]
Therefore the ideal of the degeneracy locus $Z$ is resolved
by the Eagon-Northcott complex
\[
    I_{Z/\PZ^{k+1}} 
            \from \Lambda^{k+1} {V'}^* \otimes \Lambda^{k+1} \sL
            \from \Lambda^{k+2} {V'}^* \otimes \Lambda^{k+1} \sL \otimes \sL
            \from \dots
            \from \Lambda^{2k} {V'}^* \otimes \Lambda^{k+1} \sL \otimes S_{k-1} \sL
            \from 0
\]
To show that $\iota(Z)$ is non degenerate in $\PZ^N$ we have
to prove 
\[
         h^0\bigl(I_{Z/\PZ^{k+1}}(k-2)\bigr)=0
\]
since $\iota$ is the $(k-2)$-uple embedding of $\PZ^{k+1}$.

This follows if the cohomology groups
\[
       H^i(\Lambda^{k+1} \sL \otimes S^j \sL \otimes \sO(k-2))
\]
vanish for $0 \le i \le k+1, 0 \le j < k-1$ and $0 \le i \le k, j=k-1$.
We will prove this in the next section.  
\end{proof}

\section{Cohomology of $S^i \sL (-j)$}
\nosubsections

We well calculate the cohomology of the needed homogeneous bundles
on $\PZ^{k+1}$ using the theorem of Bott. We start by fixing some notation.

Let $G=\GL(k+2)$ and $P \subset G$ the parabolic subgroup with elements
of the form
\[
\left(
\begin{array}{c|ccc}
        * & * & \dots & * \\ \hline
        0 &  \\
   \vdots &  & * \\
        0 \\
\end{array}
\right)
\]
Then $G/P \cong \PZ^{k+1}$. Let $H \subset P \subset G$ be the
subgroup of diagonal matrices, $H_i := E_{i,i}$ the natural basis
of $H$ and $\{L_i\}$ the dual basis of $H^*$. Then the $L_i$ span
the weight lattice of $G$.  The positive roots of $G$ are $L_i-L_j$ 
with $k+2 \ge i>j\ge 1$ and the fundamental weights 
are $\omega_i = \sum_{j=1}^i L_j$.

If $\rho$ is a representation of $P$, it induces a vector bundle
$E_\rho$ on $\PZ^{k+1}$ with $P$ acting on the fibers of $E_\rho$
via $\rho$. Sometimes we write 
\[
     E_\rho=E(\lambda)=E(\lambda_1, \dots, \lambda_{k+2})
\]
with $\lambda = \lambda_1 L_1 + \dots + \lambda_{k+2} L_{k+2}$
the maximal weight vector of $\rho$.

Often it is sufficient to consider the semisimple part $S_P$ of
$P$:

\begin{thm}[Classification of irreducible bundles over $G/P$]
Let $P(\Sigma) \in G$ be a parabolic subgroup and $\omega_1,\dots,\omega_k$ 
the  fundamental weights corresponding to the 
subset of simple roots $\Sigma \subset \Delta$.
Then all irreducible representations of $P(\Sigma)$ are
\[
     V \otimes L^{n_1}_{\omega_1} \otimes \dots \otimes L^{n_k}_{\omega_k} 
\]
where $V$ is a representation of $S_P$ and $n_i \in \ZZ$. $L_{\omega_i}$
are the one dimensional representations of $S_P$ induced by the 
fundamental weights.

The weight lattice of $S_P$ is embedded in the weight lattice of
$G$. If $\lambda$ is the highest weight of $V$, we will
call $\lambda + \sum n_iw_i$ the highest weight of the irreducible
representation of $P(\Sigma)$ above. 
\end{thm}

\begin{proof} \cite[Proposition 10.9 and remark 10.10]{Ott}
\end{proof}

In our case the semisimple part $S_P$ of $P$ is $\GL(1) \times  \GL(k+1)$.
Notice that the weight lattice of $\GL(1)\times\GL(k+1)$ is embedded in the
weight lattice of $\GL(k+2)$. In particular $L_1$ belongs to $\GL(1)$
and $\langle L_2,\dots,L_{k+2} \rangle$ belongs to $\GL(k+1)$.

\begin{rem}
We have $O(1)=E(1,0,\dots,0)$ since $\GL(1)$ acts on the fibers
of $O(1)$ in the standard way. In particular this representation
has maximal weight vector $L_1$. Similarity we have 
$\Omega^1(1) = E(0,1,0,\dots 0)$ since $\GL(k+1)$ acts
the fibers of $\Omega^1(1)$ with maximal weight vector $L_2$. 
Consequently
\[
       \sL = \sT_{\PZ^{k+1}}(-2) 
           = \left[ \Omega^1(2) \right]^*
           = E(1,1,0,\dots,0)^*
\]
\end{rem}

With this we are ready to use
\newcommand{\lie}[1]{{\mathfrak #1}}

\begin{thm}[Bott]
Consider the homogeneous vector bundle $E(\lambda)$ on $X=G/P$ and
$\delta$ the sum of fundamental weights of $G$. Then
\begin{itemize}
\item $H^i(X,E(\lambda))$ vanishes for all $i$ if there
is a root $\alpha$ with $(\alpha, \delta+\lambda)=0$
\item Let $i_0$ be the number of positive roots $\alpha$ with 
$(\alpha, \delta+\lambda) < 0$. Then $H^{i}(X,E(\lambda))$ vanishes
for $i\not=i_0$ and $H^{i_0}(X,E(\lambda))=\rho_{w(\delta+\lambda)-\delta}$
\end{itemize}
where $(.,.)$ denotes the Killing form on $\lie{h}^*$, 
$w(\delta+\lambda)$ is the unique element of the fundamental Weyl chamber 
which is congruent to $\delta+\lambda$ under the action of the Weyl group,
and $\rho_{w(\delta+\lambda)-\delta}$ is the corresponding representation
of $G$.
\end{thm}

\begin{proof}
\cite[Theorem 11.4]{Ott}
\end{proof}

\begin{cor}
The cohomology groups
\[
           H^i(\Lambda^{k+1} \sL \otimes S^j \sL \otimes \sO(k-2))
\]
vanish
\begin{itemize}
\item[(a)] for all $i$ if $0 \le j \le k-2$ and
\item[(b)] for all $i \not= k$ if $j=k-1$.
\end{itemize}
\end{cor}

\begin{proof}
First of all we have
\[
         \Lambda^{k+1} \sL^* = \Lambda^{k+1} \Omega^1(2) 
                           = \omega_{\PZ^{k+1}}(2k+2)
                           = \sO(k)
\]
and therefore
\[
        \Lambda^{k+1} \sL \otimes \sO(k-2) = \sO(-2) = E(-2,0,\dots,0).
\]
Similarity we get
\begin{align*}
        S^j \sL &= \left( S^j \sL^* \right)^* \\ 
              &= \left( S^j E(1,1,0,\dots,0) \right)^*\\
              &= E(j,j,0,\dots,0)^*\\
              &= E(-j,0,\dots,0,-j).
\end{align*}
Consequently
\[
   \Lambda^{k+1} \sL \otimes S^j \sL \otimes \sO(k-2) 
    = E(-j-2,0,\dots,0,-j) =: E(\lambda)
\]
Now the sum $\delta$ of the fundamental weights of $\GL(k+2)$ is
\[
       \delta = (k+2)L_1 + \dots + L_{k+2} = (k+2,k+1,\dots,2,1)
\]
so
\[
       \lambda + \delta = (k-j,k+1,\dots,2,1-j).
\]
If $0 \le j \le k-2$ then $\alpha=L_1 - L_{j+3}$ is a positive
root with $(\alpha,\lambda+\delta)=0$. Therefore, by the theorem
of Bott, all cohomology groups
vanish in this case. This proves (a).

If $j=k-1$ we have
\[ 
        \lambda + \delta = (1,k+1,\dots,2,-k)
\]
and $(\alpha,\lambda+\delta) < 0 $ for $\alpha = L_1 - L_l$ with
$2 \le l \le k+1$. All other positive roots $\alpha$ satisfy 
$(\alpha,\lambda+\delta) >0$. Consequently
\[
           H^i(\Lambda^{k+1} \sL \otimes S^j \sL \otimes \sO(k-2))
\]
vanishes for $i \not= k$ by the theorem of Bott. This proves (b).
\end{proof}

\bibliographystyle{alpha} 

\end{document}